\documentclass[10pt,a4paper,reqno]{article}
\usepackage[margin=2cm]{geometry}
\usepackage[T1]{fontenc}
\usepackage{amsmath,amssymb,amsthm,mathtools,mathrsfs}
\usepackage{charter,needspace}
\usepackage[expansion=false]{microtype}

\ifdefined\pdfgentounicode
\input{glyphtounicode}
\fi
\usepackage{xcolor,booktabs}
\usepackage[backend=biber,style=numeric,sorting=nyt,giveninits=true,
  maxbibnames=99,doi=true,url=true,isbn=false]{biblatex}
\DeclareFieldFormat[article,incollection,inproceedings,online]{title}{\mkbibemph{#1}}
\DeclareFieldFormat{journaltitle}{#1}
\DeclareFieldFormat[article]{volume}{\mkbibbold{#1}}
\renewbibmacro{in:}{\ifentrytype{article}{}{\printtext{\bibstring{in}\intitlepunct}}}
\usepackage[colorlinks=true,linkcolor=blue,citecolor=green!45!black,urlcolor=blue,bookmarksnumbered=true]{hyperref}
\hypersetup{pdftitle={Closed G2-structures with fixed metric},pdfauthor={Siva Tej Pindiprolu},pdfsubject={Closed forms with prescribed metric, compactness and invariant fibres}}
\numberwithin{equation}{section}
\allowdisplaybreaks[2]
\renewcommand{\arraystretch}{1.18}
\newcommand{\C}{\mathbb C}
\newcommand{\tr}{\operatorname{tr}}
\newcommand{\im}{\operatorname{im}}
\newcommand{\Id}{\operatorname{Id}}
\newtheorem{theorem}{Theorem}[section]
\newtheorem{mainresult}{Theorem}

\newtheorem{proposition}[theorem]{Proposition}
\newtheorem{lemma}[theorem]{Lemma}
\newtheorem{corollary}[theorem]{Corollary}
\theoremstyle{remark}\newtheorem{remark}[theorem]{Remark}
\newcommand{\R}{\mathbb R}
\newcommand{\RP}{\mathbb {RP}}
\newcommand{\Ccal}{\mathcal C}
\newcommand{\Hcal}{\mathcal H}
\newcommand{\Scal}{\operatorname{Scal}}
\newcommand{\Sym}{\operatorname{Sym}}
\newcommand{\vol}{\operatorname{vol}}
\newcommand{\av}{\mathop{\mathrm{av}}}
\newcommand{\Cl}{\operatorname{cl}}
\title{\Large Closed $G_2$-structures with fixed metric:\\compactness and almost-abelian fibres}
\author{Siva Tej Pindiprolu}
\date{\today}
\begin{document}
\setlength{\abovedisplayskip}{16pt plus 3pt minus 2pt}
\setlength{\belowdisplayskip}{16pt plus 3pt minus 2pt}
\setlength{\abovedisplayshortskip}{16pt plus 3pt minus 2pt}
\setlength{\belowdisplayshortskip}{16pt plus 3pt minus 2pt}
\maketitle
\begin{abstract}
We study the space of closed $G_2$-structures inducing a prescribed Riemannian metric and orientation. On a compact seven-manifold, with a fixed Riemannian metric and orientation, we prove that the space of compatible closed $G_2$ structures inducing a prescribed spin-structure is homeomorphic to a real projective algebraic set defined by homogeneous quadratic equations in the kernel of the corresponding Dirac operator. The defining quadrics give a local analytic obstruction map
and a \L ojasiewicz estimate that bounds the distance to the
closed locus in terms of $\|d\varphi\|_{L^2}$. We also prove smooth compactness when the prescribed metrics converge smoothly, and strong $H^1$ compactness for sequences of compatible $G_2$-structures $\varphi_{j}$ satisfying $\|d\varphi_j\|_{L^2}\to 0$. In the invariant setting, the fixed-metric fibre of any seven-dimensional almost-abelian metric Lie algebra has at most four connected components, each a compact homogeneous space. We determine these spaces spectrally when the defining endomorphism is normal, and show that for a non-normal endomorphism the fibre is empty, one circle, or two circles. On compact quotients with symmetric bracket, sufficiently small torus fibres exclude all non-invariant harmonic spinors, so the invariant classification computes the full closed fibre.
\end{abstract}

\setcounter{tocdepth}{2}
\tableofcontents

\section{Introduction and main results}
A $G_2$ structure is defined by a positive three-form on a seven-manifold which then determines a Riemannian metric and an orientation, while prescribing these data leaves a pointwise family of possible three-forms. It is therefore natural to ask which members of this family satisfy a given differential condition. We consider closedness, and study the space of positive three-forms which are closed and induce a fixed metric and orientation. The question differs from the existence of a closed $G_2$-structure on the underlying manifold, since changing the three-form is now allowed only within its metric class. It also differs from the torsion-free problem: a closed structure may have nonzero torsion, even though its defining spinor satisfies a linear elliptic equation. This last observation provides us with the starting point for the compact part of this paper.

To describe the pointwise geometry, fix the standard positive form on $\R^7$,
\[
 \varphi_{\star}=e^{127}+e^{347}+e^{567}+e^{135}-e^{146}-e^{236}-e^{245},
\]
where $e^{ijk}=e^i\wedge e^j\wedge e^k$ for the standard dual basis. A three-form on a seven-dimensional vector space is positive if it is carried to this model by a linear isomorphism. Its stabilizer is the compact Lie group $G_2$. The metric and orientation can be recovered directly from
\[
 \frac16(\iota_X\varphi)\wedge(\iota_Y\varphi)\wedge\varphi
       =g_\varphi(X,Y)\,\vol_\varphi,
\]
where $\iota_X$ denotes insertion in the first argument. In a positively oriented basis, write the left-hand side as $b_{ij}$ times the coordinate volume form. Then $g_{ij}=(\det b)^{-1/9}b_{ij}$ and $\vol_\varphi=(\det b)^{1/9}$ times that coordinate volume form. Thus the metric is determined algebraically, and not chosen independently of $\varphi$.

For a Euclidean seven-space with prescribed orientation, $SO(7)$ acts on three-forms by $R\cdot\varphi=(R^{-1})^*\varphi$. The orbit map
\[
 SO(7)/G_2\longrightarrow\{\varphi:g_\varphi=g_{\rm eucl},\ o_\varphi=o_{\rm eucl}\},
 \qquad RG_2\longmapsto(R^{-1})^*\varphi_*,
\]
is a diffeomorphism. The spin representation identifies this quotient with $\RP^7$, where $\RP^k$ is the real projective $k$-space. Bryant's same-metric parametrization \cite[equation~(3.6)]{Bryant}, also given in \cite[Theorem~2.3]{Grigorian}, makes this freedom explicit relative to a reference form. On a manifold, however, compatible forms are sections of the associated bundle with this fibre, and closedness is an equation on the section rather than on a single point of the $\RP^7$.

Throughout, $M$ denotes a smooth seven-manifold, $g$ a smooth Riemannian metric on $M$, and $o$ a chosen orientation. Let $\Omega^k(M)$ be the space of smooth real differential $k$-forms, let $d$ be the exterior derivative, and let $\Omega^3_+(M)$ be the space of positive three-forms. For $\varphi\in\Omega^3_+(M)$, write $g_\varphi$ and $o_\varphi$ for its induced metric and orientation. We distinguish the ambient fixed-metric space from its closed locus:
\[
 \mathcal P(g,o)=\{\varphi\in\Omega^3_+(M):g_\varphi=g,\ o_\varphi=o\},
 \qquad
 \Ccal(g,o)=\{\varphi\in\mathcal P(g,o):d\varphi=0\}.
\]
Both spaces carry the $C^\infty$ topology. Their points are the forms themselves, and we take no quotient by diffeomorphisms or isometries. Globally, $\mathcal P(g,o)$ is the space of smooth sections of a bundle with fibre $\RP^7$.

A $G_2$-structure is torsion-free when it is parallel for its induced Levi--Civita connection. Equivalently, $d\varphi=0$ and $d(*\varphi)=0$. Compatible torsion-free structures are described by parallel spinors, as in Grigorian \cite{Grigorian} and Lin \cite{Lin}. Their global octonion and function--one-form parametrizations describe
the sector of the reference form. To describe the entire fibre one must
also retain the other spin structures carrying parallel spinors, as we
state in Proposition \ref{prop:parallel-sectors} in Subsection \ref{sec:scalarflat}. The example
\eqref{eq:flat-two} shows why this qualification is necessary.

For a closed form with nonzero torsion, its unit defining spinor is still harmonic, but need not be parallel \cite{ACFH,BMM}. This gives us a useful finite-dimensional ambient space on a compact manifold: the Dirac kernel for the fixed metric. To recover a closed compatible form from a nonzero harmonic
spinor, we still need to require that its length be constant
and that the associated three-form be closed. Both conditions
can be expressed as quadratic equations in the spinor. We first express these constraints as a real projective algebraic set, and then use elliptic estimates to study variation of the metric and sequences whose failure of closedness tends to zero. We also determine the invariant fibre for almost-abelian metric Lie algebras: the normal case has a spectral description, while the non-normal case admits at most two circles. These results address the parameter space of forms themselves, retaining forms that become equivalent after an isometry or a Lie algebra automorphism.

\subsection{The compact fibre}
Assume here that $M$ is compact, connected and without boundary. Fix a spin structure $\sigma$, that is, a lift of the oriented orthonormal frame bundle of $(M,g,o)$ to $\operatorname{Spin}(7)$, and let $S_\sigma$ be its real rank-eight spinor bundle, with pointwise spinor inner product $\langle\ ,\ \rangle$ and norm $|s|^2=\langle s,s\rangle$. Its Levi--Civita spin connection $\nabla$ defines the Riemannian Dirac operator
\[
 D_\sigma s=\sum_{i=1}^7 e_i\cdot\nabla_{e_i}s,
 \qquad \Hcal_\sigma=\ker D_\sigma,
\]
where $(e_1,\ldots,e_7)$ is any local $g$-orthonormal frame and $\cdot$ is Clifford multiplication. Compactness makes $\Hcal_\sigma$ finite-dimensional, and elliptic regularity makes its elements smooth. The subset $\Ccal_\sigma(g,o)\subset\Ccal(g,o)$ consists of forms whose $G_2$-reduction induces the spin structure $\sigma$, and is called its spin-structure sector. A positive form induces such a lift because the inclusion $G_2\subset SO(7)$ lifts to $\operatorname{Spin}(7)$. Distinct induced spin structures are retained as distinct sectors.

For a smooth real function $f$, define its metric average by
\[
 \av(f)=\frac{1}{\vol_g(M)}\int_M f\,d\vol_g,
 \qquad \vol_g(M)=\int_M d\vol_g,
\]
where $d\vol_g$ is the Riemannian volume form for $(g,o)$. The quadratic spinor three-form is
\[
 \Phi(s)(X,Y,Z)=\varepsilon\langle\Cl(X\wedge Y\wedge Z)s,s\rangle,
\]
where $X,Y,Z$ are tangent vectors and $\Cl$ denotes the alternating Clifford product. We use the real Clifford module on which $e_1\cdots e_7$ acts as $+\Id$ for every oriented orthonormal frame, and set $\varepsilon=-1$. With this convention unit spinors yield forms inducing $(g,o)$. We use the Clifford convention $X\cdot Y+Y\cdot X=-2g(X,Y)$. The algebraic correspondence is classical \cite{Harvey,FKMS,ACFH}, and we recall the differential identities in Subsection \ref{sec:main}. Write $\mathbb P(\Hcal_\sigma)$ for the space of real one-dimensional subspaces of $\Hcal_\sigma$, and $[s]=\R s$ for the line of a nonzero spinor $s$.

\begin{mainresult}
For $(M,g,o)$ and $\sigma$ as above, the sector $\Ccal_\sigma(g,o)$ is homeomorphic to
\[
 \{[s]\in\mathbb P(\Hcal_\sigma):
       |s|^2=\av(|s|^2),\ d\Phi(s)=0\}.
\]
The homeomorphism sends $[s]$ to $\Phi(s)/\av(|s|^2)$. If $n=\dim_\R\Hcal_\sigma>0$, this is a real projective algebraic set defined by at most $n(n+1)/2$ homogeneous quadratic equations. If $n=0$, it is empty. The entire fibre is compact and has finitely many connected components.
\end{mainresult}

The theorem gives a finite-dimensional framework for the fixed-metric
problem. Its starting points are the classical correspondence between
unit spinors modulo sign and compatible forms \cite{Harvey,FKMS,ACFH},
and the fact that the defining spinor of a closed form is harmonic
\cite[Theorem 4.6]{ACFH}. Since $M$ is compact, these spinors belong
to a finite-dimensional Dirac kernel. What we use here is that the two
remaining conditions, constant spinor length and closedness of the
associated three-form, are both quadratic. This lets us describe the
whole sector by finitely many homogeneous quadratic equations and then
study its topology and local structure. The bound on the number of
equations is simply the dimension of the space of quadratic forms in
$n$ variables. The projective homeomorphism is proved in Theorem
\ref{thm:main}, and Corollary \ref{cor:all} accounts for the finitely
many spin sectors.

For the torsion-free problem, Grigorian \cite[Theorem~7.12]{Grigorian} gives a linear elliptic characterization in the reference spin sector using the Laplacian of his octonionic connection. In the closed case the Dirac equation leaves two quadratic constraints, whose Gram matrix can be formed once the harmonic spinors are known. The resulting algebra also describes projective families and small kernels, while the analytic implicit function theorem gives the local obstruction map in Subsection \ref{sec:local}.

The harmonic-spinor description also controls limits in which the prescribed metric changes. The following statement uses the $C^\infty$ topology on metrics and forms. The Sobolev space $H^1$ in its second assertion is defined using the fixed metric $g$.

\begin{mainresult}
Let $M$ be compact and connected, without boundary.
\begin{enumerate}
\item If $g_j\to g$ smoothly and $\varphi_j\in\Ccal(g_j,o)$, a subsequence converges smoothly to a member of $\Ccal(g,o)$. The map sending a closed positive form of orientation $o$ to its metric is proper, and the set of metrics admitting such a form is closed.
\item If $\varphi_j\in\mathcal P(g,o)$ and $\|d\varphi_j\|_{L^2(g)}\to0$, a subsequence converges strongly in $H^1$ to a smooth member of $\Ccal(g,o)$. In particular, if $\mathcal P(g,o)\ne\varnothing$ but $\Ccal(g,o)=\varnothing$, the closedness defect has a strictly positive infimum on $\mathcal P(g,o)$.
\end{enumerate}
\end{mainresult}

The two compactness assertions use the same unit-length normalization. For smoothly converging metrics, identifications of the spinor bundles put the harmonic-spinor equations on one fixed bundle, where uniform elliptic estimates give smooth convergence. At fixed metric, small closedness defect instead gives a small Dirac derivative in $L^2$, and Rellich compactness followed by the elliptic estimate gives strong $H^1$ convergence. The details are in Theorems \ref{thm:proper} and \ref{thm:gap}.

Combining the first assertion with the generic and local metric-deformation results of Ammann--Dahl--Humbert \cite{ADH,ADHlocal}, we obtain in Corollary \ref{cor:generic} that the set of metrics admitting a compatible closed form is closed and nowhere dense, and that compatible closed forms can be excluded by arbitrarily small metric perturbations supported in any prescribed nonempty open set.

\subsection{The invariant fibre}
Let $(V,\langle\ ,\ \rangle,o_V)$ be an oriented Euclidean six-space and $A:V\to V$ a real linear endomorphism. Set $\mathfrak g_A=V\oplus\R\nu$, give $\nu$ unit length orthogonal to $V$, and define the Lie bracket by $[\nu,X]=AX$ and $[X,Y]=0$ for $X,Y\in V$. Equip this metric Lie algebra with the orientation $o_V$ followed by $\nu$. It is almost-abelian, meaning that it contains the codimension-one abelian ideal $V$. Write $\Ccal_A$ for its positive three-forms inducing this metric and orientation and closed under the Lie algebra exterior derivative. Equivalently, these are the left-invariant compatible closed forms on the associated simply connected Lie group.

The adjoint $A^t$ is taken with respect to $\langle\ ,\ \rangle$. The endomorphism $A$ is normal if $AA^t=A^tA$ and symmetric if $A=A^t$. An orthogonal complex structure is an endomorphism $J$ satisfying $J^2=-\Id$ and $J^tJ=\Id$, where $\Id$ denotes the identity. When $[A,J]=AJ-JA=0$, the map $A$ is complex-linear on $(V,J)$, and its complex trace is denoted by $\tr_{\C,J}A$, while $\tr A$ denotes its real trace.

Closed invariant structures on nilpotent and solvable Lie groups have been studied through algebraic existence criteria and curvature restrictions, as in Conti--Fern\'andez \cite{CF} and Fern\'andez--Fino--Manero \cite{FFM}. 

In the almost-abelian setting, Freibert's existence criterion \cite[Theorem~3.1(a)]{Freibert} characterizes the Lie algebras admitting a closed form by conjugacy into $\mathfrak{sl}(3,\C)$, the complex trace-free endomorphisms of $\C^3$. We refine this criterion with the metric fixed: admissible forms correspond to orthogonal complex structures $J$ satisfying $[A,J]=0$ and $\tr_{\C,J}A=0$, together with normalized complex volume forms. Here a normalized complex volume form is a complex $(3,0)$-form $\Omega$ satisfying $\frac{i}{8}\Omega\wedge\overline{\Omega}=\vol_V$, with $i^2=-1$ and $\vol_V$ the volume form of $\langle\ ,\ \rangle$ and $o_V$, and the complex orientation of $J$ is required to equal $o_V$. The orthogonal centralizer of $A$ describes its components without a normality assumption. When $A$ is normal, the spectral decomposition makes this description explicit.

Lauret studies equivalence and solitons in this class \cite{Lauret}, while here the individual forms and the topology of their parameter space are retained. The related coclosed problem is treated by Freibert \cite{FreibertCocalibrated}, and Bagaglini--Fino \cite{BF} study its Laplacian coflow on almost-abelian groups. 

\Needspace{6\baselineskip}
\begin{mainresult}\label{main:normal}
For arbitrary $A$, the invariant fibre $\Ccal_A$ has $0$, $1$, $2$ or $4$ connected components. Each nonempty component is a projective linear subspace $\RP^1$, $\RP^3$ or $\RP^7$ of the space of invariant compatible forms, and is homogeneous under the connected orthogonal centralizer of $A$. If $A$ is not normal, its nonempty possibilities are $S^1$ and $S^1\sqcup S^1$, and both occur.

If $A$ is normal, then $\Ccal_A$ is nonempty if and only if the six complex eigenvalues of $A$, counted with algebraic multiplicity, can be written as
\[
 \zeta_1,\zeta_2,\zeta_3,
 \overline{\zeta_1},\overline{\zeta_2},\overline{\zeta_3},
 \qquad \zeta_1+\zeta_2+\zeta_3=0.
\]
This condition is independent of the chosen orientation. Each component is a principal circle bundle over a product of spaces $SO(2m)/U(m)$ and complex Grassmannians. Here $SO(2m)/U(m)$ parametrizes orthogonal complex structures of a fixed orientation on a Euclidean $2m$-space, and a complex Grassmannian $\operatorname{Gr}_p(\C^q)$ parametrizes the complex $p$-dimensional subspaces of $\C^q$. The factors and components are specified by Theorem \ref{thm:normal}, and Table \ref{tab:normal} lists every possibility.

If $A\neq0$ is normal, every nonempty component has dimension one or three. For symmetric $A$, nonemptiness is equivalent to trace zero and even-dimensional eigenspaces. For symmetric $A$, the nonempty fibre is four circles for multiplicities $2+2+2$, two copies of $\RP^3$ for multiplicities $4+2$, and $\RP^7$ for $A=0$. These multiplicities are the real dimensions of the distinct eigenspaces.
\end{mainresult}

Writing an invariant form as $\omega\wedge\nu^*+\operatorname{Re}\Omega$ reduces closedness to the condition that $A$ be complex-linear and have zero complex trace. The orthogonal centralizer then acts on the admissible complex structures, with a circle of normalized volume forms above each orbit, and a maximal torus bounds the number of components by four. The six-dimensional linear algebra in Lemma \ref{lem:two-J} gives the non-normal restriction, while the spectral decomposition determines the normal components and their stabilizers.

Theorems \ref{thm:centralizer}, \ref{thm:nonnormal} and \ref{thm:normal} give the classification, and Corollary \ref{cor:symmetric} treats the symmetric case. The standard torsion and curvature formulas of Lauret \cite{Lauret} and Moreno \cite{Moreno} relate it to the spinorial description: when $A$ is symmetric and trace-free, every invariant spinor is harmonic, although the closed locus can be a proper subset of $\RP^7$.

Lauret \cite{LauretDistinguished} describes the fixed-metric problem in the Lie-bracket formalism and asks about connectedness of the full cone of closed invariant structures on a given solvable group. Theorem \ref{main:normal} determines the fixed-metric, fixed-orientation fibre for normal almost-abelian brackets. Its disconnected examples do not decide connectedness when the metric is allowed to vary. Bazzoni--Gil-Garc\'ia \cite{BGG} instead study invariant (co)closed structures on nilmanifolds modulo Lie algebra automorphisms, without fixing the metric. That quotient is different from the fibre considered here.

The compact results are proved in Section \ref{sec:compact}, where Theorem \ref{thm:stability} gives a quantitative distance estimate from small closedness defect. Section \ref{sec:invariant} gives the invariant classification. Theorem \ref{thm:collapsed} connects the two: on compact symmetric almost-abelian quotients with sufficiently small torus fibres, it computes the whole Dirac kernel and shows that every compatible closed form is invariant.
\medskip

\noindent
\textbf{Acknowledgments.}
The author thanks Vivek Shende for his immense support,
encouragement and discussions throughout the project. The author was supported by VILLUMFONDEN, VILLUM Investigator grant 37814 and Danish National Research Foundation grant DNRF157. 

\medspace

\noindent
\textbf{AI disclosure.}
AI tools (free version of ChatGpt) were used for language editing, polishing and for literature review. The author then took the help of his friends who have access to a Pro version to assist with proof-checking. This helped improve sign and conjugation errors that were previously present, and also helped in obtaining an explicit decomposition criterion for the non-normal two-circle case, identifying the invariant components as projective linear spaces, and extending the torsion and Ricci formulas beyond the normal case, which the author has read, verified and incorporated. Rest of the ideas and proofs are from the own work of the author. The author takes full responsibility for correctness and final contents of the paper.

\section{The compact fixed-metric fibre}\label{sec:compact}
\subsection{Spinors and the quadratic constraints}\label{sec:main}
Fix a spin structure $\sigma$, and denote the corresponding sector by $\Ccal_\sigma(g,o)$. A positive three-form determines an induced spin structure: the inclusion $G_2\subset SO(7)$ lifts to $\operatorname{Spin}(7)$, since $G_2$ is simply connected. In that spin structure the form has a global unit defining spinor, unique up to a constant sign when $M$ is connected. Consequently all sectors must be retained when $H^1(M,\mathbb Z_2)$ is nonzero. A section of the projective spinor bundle in an arbitrarily fixed spin structure need not lift to a unit spinor.

More explicitly, let $F_{g,o}\to M$ be the bundle whose fibre consists of orientation-preserving linear isometries $\R^7\to T_pM$. A spin structure is a principal $\operatorname{Spin}(7)$-bundle $P_\sigma\to M$ together with a bundle map $\lambda_\sigma:P_\sigma\to F_{g,o}$ satisfying
\[
 \lambda_\sigma(pa)=\lambda_\sigma(p)\rho(a),
 \qquad \rho:\operatorname{Spin}(7)\longrightarrow SO(7),
\]
where $\rho$ is the double covering. The spinor bundle is $P_\sigma\times_{\operatorname{Spin}(7)}\Delta_7$, where $\Delta_7$ is the real eight-dimensional spin representation. If $P_\varphi\subset F_{g,o}$ is the bundle of frames sending the model form to $\varphi$, then $P_\varphi$ is a principal $G_2$-bundle. The lifted inclusion $G_2\hookrightarrow\operatorname{Spin}(7)$ constructs its induced spin bundle as $P_\varphi\times_{G_2}\operatorname{Spin}(7)$. The $G_2$-fixed unit vector in $\Delta_7$ then gives the defining unit spinor. Conversely, its stabilizer reconstructs $P_\varphi$. This is the standard spinorial description of a $G_2$-reduction \cite{Gray,Lawson,CN}. It explains why the spin structure carried by the form must be used when passing from a projective spinor to a global unit spinor.

Let $S_\sigma$ be the real rank-eight spinor bundle, $\nabla$ its Levi--Civita spin connection, and
\[
 D_\sigma=\sum_{i=1}^7 e_i\mathbin{\cdot}\nabla_{e_i},\qquad
 \Hcal_\sigma=\ker D_\sigma,\qquad n=\dim_\R\Hcal_\sigma.
\]
Here $X\cdot Y+Y\cdot X=-2g(X,Y)$. Assume throughout this section that $M$ is compact, connected and without boundary. Elliptic theory makes this kernel finite-dimensional and smooth. All spinor dimensions below are real.

Using the metric average defined in the introduction, set
\[
 (s,t)_0=\av\langle s,t\rangle.
\]
Recall the quadratic three-form $\Phi(s)$:
\begin{equation}
 \Phi(s)(X,Y,Z)=\varepsilon\langle\Cl(X\wedge Y\wedge Z)s,s\rangle.
\end{equation}
We retain the convention $\varepsilon=-1$ fixed in the introduction. For a unit spinor $u$, $\Phi(u)$ induces exactly $g$ and $o$. The fibrewise correspondence is $u\bmod\{\pm1\}\leftrightarrow\Phi(u)$ \cite{Harvey,FKMS,ACFH,Grigorian}. Its polarization is
\[
 B(s,t)=\tfrac12\bigl(\Phi(s+t)-\Phi(s)-\Phi(t)\bigr),\qquad B(s,s)=\Phi(s).
\]
Both $B$ and the spinor inner product are symmetric real bilinear maps.

For a unit spinor $u$, put $\varphi=\Phi(u)$ and $\psi=*_{g,o}\varphi$. We normalize the full torsion tensor by
\[
 \nabla_i\varphi_{jkl}=T_i{}^m\psi_{mjkl},\qquad
 T(X)=\bigl(T(X,\cdot)\bigr)^\sharp.
\]
This is the full-torsion normalization of \cite[Theorem~2.27]{Karigiannis}, with the Hodge convention $*(\beta\wedge\varphi)=-\beta$ on $\Lambda^2_{14}$ used here. The relation to the spin connection can be checked directly. Clifford multiplication identifies $T_pM$ isometrically with $u(p)^\perp$, and the algebraic spinor-to-form identity is
\begin{equation}\label{eq:variation-spin}
 B(u,Y\cdot u)=-\iota_Y\psi.
\end{equation}
For example, this follows by evaluating the alternating Clifford product in an adapted orthonormal frame. Since $|u|=1$, there is a unique vector $a_X$ such that $\nabla_Xu=a_X\cdot u$. Differentiating $\varphi=\Phi(u)$ now gives
\[
 \nabla_X\varphi=2B(u,a_X\cdot u)=-2\iota_{a_X}\psi.
\]
Comparison with the definition of $T$ proves
\[
 \nabla_Xu=-\tfrac12T(X)\cdot u.
\]
Here $T(X)$ is a tangent vector, and the dot denotes Clifford multiplication. We will use this relation throughout, with the factor and sign determined by the conventions fixed above.

We use the torsion forms defined by
\[
 d\varphi=\tau_0\psi+3\tau_1\wedge\varphi+*\tau_3,
 \qquad d\psi=4\tau_1\wedge\psi+\tau_2\wedge\varphi,
\]
where $\tau_0$ is a function, $\tau_1$ a one-form, $\tau_2\in\Omega^2_{14}$ and $\tau_3\in\Omega^3_{27}$. Therefore closedness is equivalent to $\tau_0=\tau_1=\tau_3=0$, or to $T\in\Lambda^2_{14}\cong\mathfrak g_2$, in which case $T_{ij}=-\frac12(\tau_2)_{ij}$ in the conventions above \cite{LW}. The decomposition into torsion classes is the Fern\'andez--Gray classification \cite{FG}.

Introduce the symmetric bilinear constraint
\begin{equation}
 Q(s,t)=\bigl(\langle s,t\rangle-(s,t)_0,\ dB(s,t)\bigr)
 \in C^\infty(M)\oplus\Omega^4(M).
\end{equation}
Equivalently, let $\mathscr Q:\Sym^2\Hcal_\sigma\to C^\infty(M)\oplus\Omega^4(M)$ be the linear map with $\mathscr Q(s\odot t)=Q(s,t)$, where $s\odot t=(s\otimes t+t\otimes s)/2$.

\begin{theorem}\label{thm:main}
Let $(M^7,g,o)$ be compact and connected, without boundary, and fix a spin structure $\sigma$. Then
\begin{equation}\label{eq:Z}
 Z_\sigma=\{[s]\in\mathbb P(\Hcal_\sigma):Q(s,s)=0\}
 \ \longrightarrow\ \Ccal_\sigma(g,o),\qquad
 [s]\longmapsto\frac{\Phi(s)}{(s,s)_0}
\end{equation}
is a homeomorphism, where the right-hand side has its $C^\infty$ topology. If $n=0$, the solution set is empty.If $n>0$, then, after choosing a basis of $\Hcal_\sigma$,
the set $Z_\sigma\subset\RP^{n-1}$ is defined by at most
$n(n+1)/2$ homogeneous real quadratic equations.

In particular, $\Ccal_\sigma(g,o)$ is compact, has finitely many connected components, and has real semialgebraic dimension at most $n-1$. No real analyticity assumption on $g$ is required.
\end{theorem}

\begin{proof}
We first recall why a closed structure has a harmonic defining spinor, following the spinorial torsion calculation of Agricola--Chiossi--Friedrich--H\"oll \cite{ACFH}. In an orthonormal frame,
\[
 D_\sigma u=\tfrac12\tr(T)u-
 \tfrac12\sum_{i<j}(T_{ij}-T_{ji})e_i\cdot e_j\cdot u.
\]
The symmetric part contributes its trace. The skew part acts through its $\Lambda^2_7$ component because $\Lambda^2_{14}=\mathfrak g_2$ annihilates $u$.
Consequently a closed form has $D_\sigma u=0$. More precisely,
\[
 D_\sigma u=0\quad\Longleftrightarrow\quad\tau_0=\tau_1=0,
\]
as in \cite[Theorem~4.6]{ACFH} and \cite[Corollary~9.5]{Grigorian}. Among these unit harmonic spinors, closedness additionally requires $\tau_3=0$, so a nonzero $\tau_3$ gives a non-closed structure. Proposition \ref{prop:dirac-aa} supplies examples.

If $[s]\in Z_\sigma$, put $c^2=(s,s)_0>0$. The first component of $Q(s,s)=0$ says $|s|^2=c^2$ everywhere. Hence $u=s/c$ is a global unit spinor. The second component says $d\Phi(s)=0$, so $d\Phi(u)=0$. The form $\Phi(u)$ has metric $g$, orientation $o$, and induced spin structure $\sigma$. This proves that \eqref{eq:Z} is well defined, since scaling $s$ by any nonzero real constant has no effect.

Conversely, a form in $\Ccal_\sigma(g,o)$ has a defining unit spinor $u$. The preceding calculation puts $u$ in $\Hcal_\sigma$, and both components of $Q(u,u)$ vanish. If two normalized spinors define the same form, they differ pointswise by sign. Their nonvanishing and connectedness of $M$ make this sign constant. Thus \eqref{eq:Z} is bijectve.

The domain of $\mathscr Q$ has dimension $d_{\mathrm{sym}}=n(n+1)/2$, so its image is a vector space of dimension $r\leq d_{\mathrm{sym}}$, even though its target is infinite-dimensional. Choose a basis $F_1,\ldots,F_r$ of this image. In any basis of $\Hcal_\sigma$ one can write
\[
 Q(s,s)=\sum_{\mu=1}^r q_\mu(s)F_\mu
\]
with homogeneous quadratic polynomials $q_\mu$. Linear independence now lets us replace the quadratic constraint with the claimed finite collection of scalar equations. The resulting zero set is closed in a finite-dimensional real projective space, and is therefore compact. It remains to check that the bijection is a homeomorphism.
The expression in \eqref{eq:Z} involves only finitely many fixed
smooth spinors, and its denominator is a positive definite
quadratic form. It therefore depends continuously on the
projective spinor in the $C^\infty$ topology. Since the domain
is compact and the space of smooth three-forms is Hausdorff,
this continuous bijection is a homeomorphism. Finally, a real projective algebraic set is semialgebraic and has finitely many connected components, and its dimension is bounded by that of $\RP^{n-1}$, by the standard finiteness properties of semialgebraic sets \cite{BCR}.
\end{proof}

\begin{corollary}\label{cor:all}
For fixed $(g,o)$ on compact connected $M$, the whole space $\Ccal(g,o)$ is a finite disjoint union of the algebraic models $Z_\sigma$ above. In particular it is compact and has finitely many connected components. If $M$ is not spin, it is empty.
\end{corollary}
\begin{proof}
Every compatible $G_2$-structure has its uniquely determined induced spin structure. When spin structures exist, their equivalence classes form a torsor for the finite group $H^1(M,\mathbb Z_2)$. Therefore, there are finitely many sectors, and the maps above give a continuous bijection from their compact disjoint union to $\Ccal(g,o)$. The compact-to-Hausdorff argument at the end of Theorem \ref{thm:main} applies to this disjoint union.
\end{proof}

\subsection{Quadratic equations and the local structure}\label{sec:local}
Choose a basis $s_1,\ldots,s_n$ orthonormal for $(\ ,\ )_0$. For indices $i\leq j$, set
\[
 F_{ij}=Q(s_i,s_j)=
 \bigl(\langle s_i,s_j\rangle-\delta_{ij},\ dB(s_i,s_j)\bigr).
\]
Equip $C^\infty(M)\oplus\Omega^4(M)$ with the positive definite inner product
\[
 ((f,\eta),(h,\xi))_{\mathrm{av}}=\av\!\left(fh+\langle\eta,\xi\rangle_g\right).
\]
The $d_{\mathrm{sym}}\times d_{\mathrm{sym}}$ real Gram matrix and the degreetwo coordinate vector are
\begin{equation}
 G_{ij,k\ell}=(F_{ij},F_{k\ell})_{\mathrm{av}},\qquad
 v(a)_{ii}=a_i^2,\quad v(a)_{ij}=2a_i a_j\quad(i<j).
\end{equation}

\begin{proposition}
For $a\in\R^n\setminus\{0\}$, the following conditions are equivalent:
\[
 \left[\sum_i a_i s_i\right]\in Z_\sigma,
 \qquad Gv(a)=0,
 \qquad v(a)^t Gv(a)=0.
\]

When they hold, the corresponding closed form is
\begin{equation}\label{eq:reconstruct}
 \varphi_a=
 \frac{\displaystyle\sum_i a_i^2\Phi(s_i)
       +2\sum_{i<j}a_i a_j B(s_i,s_j)}{\displaystyle\sum_i a_i^2}.
\end{equation}
Hence nonemptiness is exactly the existence of a real nonzero $a$ satisfying $Gv(a)=0$. Equivalently, $\ker\mathscr Q$ contains a nonzero rank-one positive semidefinite tensor $s\otimes s$.
\end{proposition}
\begin{proof}
If $s=\sum_i a_i s_i$, then $Q(s,s)=\sum_{i\leq j}v(a)_{ij}F_{ij}$. Consequently
\[
 v(a)^tGv(a)=\|Q(s,s)\|_{\mathrm{av}}^2.
\]
A smooth function or form of zero $L^2$ norm vanishes identically. Since $G$ is positive semidefinite, its quadratic form vanishes precisely on its kernel. Theorem \ref{thm:main} proves the equivalences, and polarization gives \eqref{eq:reconstruct}.
\end{proof}

This is a finite algebraic formulation once a basis of $\ker D_\sigma$ and its Gram matrix have been determined. Positive definiteness of $G$ excludes solutions, while a singular matrix gives a solution only when its kernel meets the rank-one cone described above.

Once the harmonic spinors have been fixed, the remaining constraints are homogeneous quadrics. In small kernel dimension their common zero set is particularly restricted. The following consequences record what follows from this algebra alone, without additional assumptions on the torsion of a chosen solution.

\Needspace{5\baselineskip}
\begin{corollary}
In one fixed spin structure:
\begin{enumerate}
 \item If $n=0$, there are no compatible closed $G_2$-structures.
 \item If $n=1$, there is at most one. For a nonzero generator $s$, it exists exactly when $|s|$ is constant and $d\Phi(s)=0$.
 \item If $n=2$, the solution set is empty, one point, two points, or all of $\RP^1$. In particular, three distinct solutions imply an entire projective circle of solutions.
\end{enumerate}
In the case $n=2$, the entire circle occurs if and only if $G=0$.
\end{corollary}
\begin{proof}
The first two assertions follow from Theorem \ref{thm:main}. For $n=2$, if all the defining quadratic polynomials vanish identically, the zero locus is $\RP^1$. Otherwise it is contained in the zeros of a nonzero homogeneous quadratic in two variables, which has at most two distinct real projective zeros. All quadratic constraints vanish identically if and only if $Q(s,s)=0$ for every $s\in\Hcal_\sigma$. Polarization makes this equivalent to $Q(s,t)=0$ for every pair, or $F_{ij}=0$ for all $i,j$, which is equivalent to $G=0$.
\end{proof}

\begin{proposition}
Let $W\subset\Hcal_\sigma$ have dimension $k>0$. Every line in $W$ gives a compatible closed form if and only if
\begin{equation}\label{eq:family}
 \langle s,t\rangle\ \text{is constant on }M,
 \qquad dB(s,t)=0\qquad(s,t\in W).
\end{equation}
It suffices to check these conditions on all pairs of a basis of $W$. In this case the resulting subset is homeomorphic to $\RP^{k-1}$, and necessarily $k\leq8$.
\end{proposition}
\begin{proof}
The condition that every line is a solution is $Q(s,s)=0$ for every $s\in W$. Polarization makes this equivalent to $Q|_{W\times W}=0$, which is precisely \eqref{eq:family}. Restricting the homeomorphism of Theorem \ref{thm:main} to
$\mathbb P(W)$ gives the claimed projective family. To obtain
the dimension bound, fix a point $p\in M$ and consider the
evaluation map $W\to(S_\sigma)_p$. If a spinor in $W$ vanishes
at $p$, its constant length must be zero, so it vanishes
everywhere. The evaluation map is therefore injective, and
since the spinor fibre has real dimension eight, we obtain
$k\leq8$.
\end{proof}

To describe the fibre near a particular closed form, choose its defining harmonic spinor $u\in\Hcal_\sigma$, normalized by $(u,u)_0=1$. Since it is a solution, this normalization also makes its pointwise length one. We now choose an affine chart around $[u]$ and study the constraint in these coordinates. The first step is to solve its component in the image of the linearization. The part that remains will give us the obstruction map. Set
\[
 \mathcal X_u=u^\perp\cap\Hcal_\sigma,\qquad \mathcal Y=\im\mathscr Q,
 \qquad \ell_u:\mathcal X_u\longrightarrow \mathcal Y,\quad \ell_u b=2Q(u,b).
\]
Give $\mathcal Y$ the inner product $(\ ,\ )_{\mathrm{av}}$ defined above, and write $\|s\|_0=(s,s)_0^{1/2}$. In the affine projective chart $[u+b]$, $b\in \mathcal X_u$, the exact equation is
\begin{equation}\label{eq:affine}
 F(b)=\ell_u b+Q(b,b)=0.
\end{equation}
No terms of degree greater than two occur in this chart.
Write $\mathcal K=\ker \ell_u$, $\mathcal N=\mathcal K^\perp\subset \mathcal X_u$, $\mathcal I=\im \ell_u$, and $\mathcal O=\mathcal I^\perp\subset \mathcal Y$.
Let $P_{\mathcal I},P_{\mathcal O}$ be the orthogonal projections and $\mathcal R_u:\mathcal I\to \mathcal N$ the inverse of $\ell_u|_{\mathcal N}$.

\begin{theorem}\label{thm:local}
There is a real analytic map $h$ from a neighbourhood of $0$
in $\mathcal K$ to $\mathcal N$, unique as a germ at $0$,
with $h(0)=0$, such that
\[
 P_{\mathcal I}F(k+h(k))=0.
\]
The germ of $\Ccal_\sigma(g,o)$ at $\Phi(u)$ is homeomorphic to the zero germ of
\begin{equation}\label{eq:kappa}
 \kappa:\mathcal K\longrightarrow \mathcal O,\qquad
 \kappa(k)=P_{\mathcal O} Q(k+h(k),k+h(k)).
\end{equation}
Its first terms are
\begin{align}
 h(k)&=-\mathcal R_u P_{\mathcal I}Q(k,k)+O(\|k\|_0^3),\label{eq:h-exp}\\
 \kappa(k)&=P_{\mathcal O}Q(k,k)-2P_{\mathcal O}Q(k,\mathcal R_u P_{\mathcal I}Q(k,k))+O(\|k\|_0^4).
 \label{eq:k-exp}
\end{align}
In particular, $\mathcal O=0$ implies that the solution germ is a smooth real analytic manifold of dimension $\dim \mathcal K$.
\end{theorem}
\begin{proof}
This is the finite-dimensional Lyapunov--Schmidt reduction of the quadratic equation \eqref{eq:affine}. We give the expansion because it will distinguish the tangent equation from its second-order obstruction. Apply the real analytic implicit function theorem to
\[
 (k,n)\longmapsto \ell_u n+P_{\mathcal I}Q(k+n,k+n),\qquad \mathcal K\times \mathcal N\to \mathcal I.
\]
Its derivative in $n$ at $(0,0)$ is the isomorphism $\ell_u|_{\mathcal N}$. This gives the unique germ $h$ with $h(0)=0$. The derivative in $k$ of the range equation vanishes at $(0,0)$, so differentiation gives $dh_0=0$. Thus $h=O(\|k\|^2)$. Substitution into the range equation yields
\[
 h(k)=-\mathcal R_u P_{\mathcal I}Q(k+h(k),k+h(k)),
\]
which proves \eqref{eq:h-exp}. Every sufficiently small solution of the range equation is uniquely of this form. Its remaining equation is exactly \eqref{eq:kappa}, since $P_{\mathcal O}\ell_u=0$. Expanding the bilinear expression and using \eqref{eq:h-exp} proves \eqref{eq:k-exp}, since the term $Q(h,h)$ is of order four. If $\mathcal O=0$, there is no remaining equation and the solution set is the analytic graph of $h$.
\end{proof}

\begin{corollary}
A direction $k\in \mathcal K$ extends to a second-order solution
$b(t)=tk+t^2b_2+O(t^3)$ with $F(b(t))=O(t^3)$ if and only if
\[
 P_{\mathcal O}Q(k,k)=0.
\]
If $P_{\mathcal O}Q(k,k)\ne0$ for every nonzero $k\in \mathcal K$, then $\Phi(u)$ is isolated, even when $\mathcal K$ is nonzero.
\end{corollary}
\begin{proof}
The coefficient of $t^2$ in $F(b(t))$ is $\ell_u b_2+Q(k,k)$. It vanishes precisely when $Q(k,k)\in \mathcal I$, in which case $b_2=-\mathcal R_u P_{\mathcal I}Q(k,k)$ is a choice. This gives the required second-order solution. To prove the
isolation statement, first suppose that $\mathcal K\ne0$.
By hypothesis, $P_{\mathcal O}Q(k,k)$ is nonzero on the unit
sphere of $\mathcal K$. Compactness of that sphere, together
with the quadratic homogeneity, gives a constant $c>0$ such that
\[
 \|P_{\mathcal O}Q(k,k)\|_{\mathrm{av}}\ge c\|k\|_0^2.
\]
The remaining terms in \eqref{eq:k-exp} are of higher order,
so they cannot cancel this quadratic term when $k$ is
sufficiently small and nonzero. Hence $\kappa$ has no other
zero near the origin, and $\Phi(u)$ is isolated. If
$\mathcal K=0$, the local model already consists of a single
point.
\end{proof}

Karigiannis \cite{KarDeform} describes fixed-metric variations of
$G_2$-structures in terms of vector fields. Using this description,
we can express the linearized closedness condition as an equation
for a vector field. This gives a geometric interpretation of the
kernel $\mathcal K$ in Theorem \ref{thm:local}.

\begin{proposition}
Let $u$ define $\varphi\in\Ccal_\sigma(g,o)$, and put $\psi=*\varphi$. Under $X\mapsto X\cdot u$, the kernel $\mathcal K$ in Theorem \ref{thm:local} is
\[
 \mathcal K\cong\{X\in\Gamma(TM):d(\iota_X\psi)=0\}
 =\{X:\mathcal L_X\psi=\iota_X(\tau_2\wedge\varphi)\}.
\]
More generally, writing $b=fu+X\cdot u$, the linearized constraint is
\begin{equation}\label{eq:geometric-linear}
 \ell_u b=\bigl(2f,\,2d(f\varphi-\iota_X\psi)\bigr),
 \qquad D_\sigma b=0,\quad\av(f)=0.
\end{equation}
Consequently $\mathcal O$ is naturally isomorphic to $\mathcal Y$ modulo the image of the map in \eqref{eq:geometric-linear}.
\end{proposition}
\begin{proof}
The first component of $\ell_u b=0$ makes $b$ pointwise perpendicular to $u$. Therefore $b=X\cdot u$ for a unique vector field, and \eqref{eq:variation-spin} identifies its second component with $-2d(\iota_X\psi)$. Conversely, suppose this exterior derivative vanishes. The unit spinor path
\[
 u_t=\frac{u+tX\cdot u}{(1+t^2|X|^2)^{1/2}}
\]
has $d\Phi(u_t)=O(t^2)$. Lemma \ref{lem:dirac-bound}, proved below by the pointwise torsion decomposition, gives $D_\sigma u_t=O(t^2)$. Differentiating at zero shows $D_\sigma(X\cdot u)=0$, which is the harmonicity condition required in the definition of $\mathcal K$. Cartan's formula and $d\psi=\tau_2\wedge\varphi$ give the second description. For general $b$, polarzation gives $B(u,b)=f\varphi-\iota_X\psi$, proving \eqref{eq:geometric-linear}. Finally, the orthogonal complement of $\im\ell_u$ in $\mathcal Y$ identifies with its quoteint, as in the reduction above.
\end{proof}

When $\varphi$ is torsion-free, this kernel consists of the infinitesimal automorphisms of $\psi$. Such a vector fiel preserves $g$, since $\psi$ determines the metric, and is therefore Killing. On a compact Ricci-flat manifold a Killing vector field is parallel by the Bochner identity. This recovers the tangent space of the projective family in \eqref{eq:lin}. For nonzero torsion the equation instead has the explicit right-hand side displayed above. It does not by itself establish the existence of a singular or obstructed fiber.

\begin{proposition}
Two forms in $\Ccal(g,o)$ which agree on a nonempty open set agree everywhere on connected $M$.
\end{proposition}
\begin{proof}
Lift the forms to the universal cover. Their induced spin structures are isomorphic there, and their defining unit spinors can be chosen to agree over one lift of the given open set. The difference is a harmonic spinor. Weak unique continuation for the Dirac operator \cite{BarNodal} makes this difference zero on the connected cover. The forms therefore agree on $M$.
\end{proof}

This argument applies across spin sectors, because the comparison takes place on the universal cover. It proves uniqueness from an open set, not from the value at a single point. Injectivity of point evaluation for a general closed fibre, and a corresponding universal dimension bound of seven, remain separate questions.

\subsection{Compactness and the closedness defect}
The finite-dimensional model describes a single fibre, but its construction uses the kernel of a Dirac operator whose dimension can jump when the metric changes. It is therefore useful to formulate a compactness statement which does not require these kernels to form a vector bundle. The point is that a closed compatible form supplies a harmonic spinor of pointwise unit length. This normalization gives the uniform bound needed for elliptic estimates, even at a metric where the dimension of the kernel increases.

Write $\operatorname{Met}^\infty(M)$ for the space of smooth Riemannian metrics with its $C^\infty$ topology, and put
\[
 \mathscr Z_o=\{\varphi\in\Omega^3_+(M):o_\varphi=o,\ d\varphi=0\},
 \qquad \mathfrak m:\mathscr Z_o\longrightarrow\operatorname{Met}^\infty(M),
 \quad \mathfrak m(\varphi)=g_\varphi.
\]
Thus $\mathfrak m^{-1}(g)=\Ccal(g,o)$. Both spaces are metrizable with these topologies, so compactness can be tested using sequences.

\begin{theorem}\label{thm:proper}
Let $M$ be compact, connected and without boundary. Suppose $g_j\to g$ in $C^\infty$ and $\varphi_j\in\Ccal(g_j,o)$. Then a subsequence converges in $C^\infty$ to a form $\varphi\in\Ccal(g,o)$. Consequently $\mathfrak m$ is proper, in the sense that the inverse image of every compact subset of $\operatorname{Met}^\infty(M)$ is compact, and its image is closed.
\end{theorem}
\begin{proof}
There are finitely many spin structures on $M$, so after passing to a subsequence the forms induce the same spin structure $\sigma$. Choose their unit defining spinors $u_j$. The closed-to-harmonic implication recalled in the proof of Theorem \ref{thm:main} gives $D_{g_j,\sigma}u_j=0$. The metric identifications of spinor bundles introduced by Bourguignon–Gauduchon allow us to regard these equations as equations on a single bundle \cite{BGmetric}. More explicitly, there is a unique positive $g$-self-adjoint bundle map $b_j:TM\to TM$ such that
\[
 g_j(b_jX,b_jY)=g(X,Y).
\]
It sends an oriented $g$-orthonormal frame to a $g_j$-orthonormal frame. Its lift to the spin structures gives a fibrewise isometry $\beta_j:S_{g,\sigma}\to S_{g_j,\sigma}$. Either choice of the overall sign of the lift can be used. The maps $b_j$ converge smoothly to the identity, and on the fixed bundle $S_{g,\sigma}$ the operators
\[
 \widehat D_j=\beta_j^{-1}D_{g_j,\sigma}\beta_j
\]
converge in all coefficient derivatives to $D_{g,\sigma}$, as in the metric comparison in \cite[Section~2.1]{ADH}.

Set $v_j=\beta_j^{-1}u_j$. Then $|v_j|_g=1$ and $\widehat D_jv_j=0$. Uniform elliptic estimates on the compact manifold give, for each integer $r\geq0$ and all sufficiently large $j$,
\[
 \|v_j\|_{H^{r+1}(g)}\leq C_r\bigl(\|\widehat D_jv_j\|_{H^r(g)}+
                                      \|v_j\|_{L^2(g)}\bigr).
\]
The constants can be chosen uniformly because the coefficients
converge smoothly and the limiting principal symbol is elliptic.
Since $|v_j|_g=1$, we also have
$\|v_j\|_{L^2(g)}=\vol_g(M)^{1/2}$. The estimates therefore give
a uniform bound in every Sobolev norm. Applying Sobolev
compactness and taking a diagonal subsequence, we obtain
smooth convergence $v_j\to v$. Passing to the limit preserves
the unit-length condition and gives $D_{g,\sigma}v=0$.

The spinor-to-form construction is natural under the metric identifications, and is quadratic in the spinor. Hence $\varphi_j\to\Phi_g(v)$ smoothly. The limit has metric $g$ and orientation $o$, and it is closed because exterior differentiation is continuous in this topology. This proves the sequential assertion. If the metrics range in a compact set, first choose a convergent subsequence of metrics and then apply this assertion, which proves properness. Applying it directly to a convergent sequence in the image proves that the image is closed.
\end{proof}

There is also a finite-regularity version. If $g_j\to g$ in $C^1$, with $g$ smooth and $\varphi_j\in\Ccal(g_j,o)$, a subsequence converges strongly in $H^1(g)$ to a smooth member of $\Ccal(g,o)$. To see this, use the same spinor identifications as in the proof. They give
\[
 \|\widehat D_j-D_{g,\sigma}\|_{H^1\to L^2}\longrightarrow0.
\]
The elliptic estimate for the fixed operator absorbs this small difference and gives a uniform $H^1$ bound on $v_j$. Rellich compactness gives a unit weak limit $v$ with $D_{g,\sigma}v=0$. Moreover,
\[
 D_{g,\sigma}(v_j-v)=(D_{g,\sigma}-\widehat D_j)v_j\longrightarrow0
 \quad\hbox{in }L^2,
\]
so the same estimate gives strong $H^1$ convergence. The bilinear differentiation argument in Theorem \ref{thm:gap} then applies to the forms. The metric identifications converge in $C^1$, so they do not change this conclusion. Closedness passes to the distributional limit, and elliptic regularity makes the limiting spinor smooth. The smoothness of $g$ is used here, even though convergence of the metrics is only in $C^1$.

These statements assume convergence of the metrics as tensors on a fixed manifold. Lotay--Wei \cite{LW} prove a geometric compactness theorem for Laplacian-flow solutions modulo diffeomorphisms under curvature and injectivity-radius hypotheses. Their theorem obtains the geometric convergence that is assumed here.

\begin{corollary}\label{cor:generic}
The subset
\[
 \mathscr A_o=\{g\in\operatorname{Met}^\infty(M):\Ccal(g,o)\ne\varnothing\}
\]
is closed and nowhere dense in every $C^k$ topology on smooth metrics, $1\leq k\leq\infty$. Moreover, for every metric $g$, every nonempty open subset $U\subset M$, and every $C^\infty$ neighbourhood of $g$, there is a metric $g'$ in that neighbourhood, equal to $g$ outside $U$, with $\Ccal(g',o)=\varnothing$.
\end{corollary}
\begin{proof}
The non-spin case is immediate. Otherwise, in dimension seven the index-theoretic lower bound for the Dirac kernel is zero. The generic-metric theorem of Ammann--Dahl--Humbert \cite[Theorem~1.1]{ADH} therefore gives an open dense set of metrics with invertible Dirac operator in each spin structure. There are only finitely many spin structures, so the intersection of these sets is still open and dense. By Theorem \ref{thm:main}, its metrics admit no compatible closed form. The finite-regularity argument following Theorem \ref{thm:proper} supplies closedness in $C^1$, and hence in each of these topologies. The genericity sets are $C^1$-open and dense in each finite $C^k$ topology, and also dense in $C^\infty$. This proves nowhere density. The dependence of harmonic spinors on the metric goes back to Hitchin \cite{HitchinSpin} and B\"ar \cite{BarMetrics}.

For the local statement, choose a ball with closure in $U$. The local deformation theorem \cite[Corollary~1.2]{ADHlocal} gives openness in $C^1$ and density in every finite $C^k$ topology, $k\geq1$, among metrics equal to $g$ outside that ball. Intersect these open dense sets for the finitely many spin structures. A $C^\infty$ neighbourhood contains a neighbourhood specified by finitely many derivative bounds, so the same finite-intersection argument supplies $g'$ in the prescribed neighbourhood.
\end{proof}

Hence, existence is closed under smooth limits of metrics, although a nonempty fibre can disappear under arbitrarily small perturbations, including perturbations supported in a prescribed open set.

There is a related statement in which the metric is fixed and closedness is approached only in $L^2$. For $\varphi\in\mathcal P(g,o)$ call $\|d\varphi\|_{L^2(g)}$ its closedness defect. Pointwise, $d\varphi$ does not detect the $\Lambda^2_{14}$ torsion
component. On a compact manifold with fixed metric, however, its
$L^2$ norm controls the full torsion energy through the integrated
identity \eqref{eq:torsion-energy} below. For the compactness argument,
we first use its control of the Dirac operator of the defining unit
spinor.

\begin{lemma}\label{lem:dirac-bound}
Every smooth unit spinor $u$ satisfies
\[
 |D_\sigma u|^2=\frac{49}{64}\tau_0^2+9|\tau_1|^2,
 \qquad |D_\sigma u|\leq\tfrac12|d\Phi(u)|.
\]
Here the torsion forms are those of $\Phi(u)$.
\end{lemma}
\begin{proof}
The scalar and vector components of the Dirac formula in Theorem \ref{thm:main} have lengths $7|\tau_0|/8$ and $3|\tau_1|$, respectively. This is the usual spinorial torsion formula \cite{ACFH,Grigorian}, in the normalization of the exterior equations above. The components are orthogonal since $u\perp X\cdot u$. On the other hand, orthogonality of the irreducible four-form components gives
\[
 |d\Phi(u)|^2=7\tau_0^2+36|\tau_1|^2+|\tau_3|^2.
\]
Comparing the coefficients in these two expressions gives
$|D_\sigma u|^2\leq\frac14|d\Phi(u)|^2$, which is the desired
inequality. For torsion with only a nonzero $\tau_1$ component,
equality holds, so the constant $1/2$ is sharp as a pointwise
algebraic bound.
\end{proof}

If $D_\sigma$ is invertible and its least absolute eigenvalue is $\mu_\sigma>0$, integration gives the explicit bound
\begin{equation}
 \|d\varphi\|_{L^2}\geq2\mu_\sigma\vol_g(M)^{1/2}
 \qquad(\varphi\in\mathcal P_\sigma(g,o)).
\end{equation}
We write $\mathcal P_\sigma(g,o)$ for compatible forms inducing $\sigma$. This estimate uses the full Dirac spectrum, while the qualitative gap below also applies when the kernel is nonzero.

\begin{theorem}\label{thm:gap}
Let $(M,g,o)$ be compact and connected, without boundary. If $\varphi_j\in\mathcal P(g,o)$ and $\|d\varphi_j\|_{L^2(g)}\to0$, then a subsequence converges strongly in $H^1$ to a smooth form in $\Ccal(g,o)$. In particular, if $\mathcal P(g,o)$ is nonempty and $\Ccal(g,o)$ is empty, then
\[
 \inf_{\varphi\in\mathcal P(g,o)}\|d\varphi\|_{L^2(g)}>0.
\]
Here $H^1$ is the Sobolev space defined using the fixed metric and its Levi--Civita connection.
\end{theorem}
\begin{proof}
As in the proof of Theorem \ref{thm:proper}, we first pass to a
subsequence for which all the forms induce the same spin structure
$\sigma$, and choose unit spinors $u_j$ defining $\varphi_j$.
Lemma \ref{lem:dirac-bound} then gives
$D_\sigma u_j\to0$ in $L^2$. The elliptic estimate
\[
 \|u_j\|_{H^1}\leq C\bigl(\|D_\sigma u_j\|_{L^2}
                         +\|u_j\|_{L^2}\bigr),
\]
together with $|u_j|=1$, gives a uniform $H^1$ bound. Rellich compactness gives a subsequence converging weakly in $H^1$, strongly in $L^2$, and almost everywhere to a spinor $u$. Therefore, $|u|=1$ almost everywhere and $D_\sigma u=0$ distributionally. Elliptic regularity makes $u$ smooth and hence unit everywhere.

The quadratic map $\Phi$ is uniformly Lipschitz on the unit sphere in each spinor fibre. Consequently $\varphi_j\to\Phi(u)$ in $L^2$. Although we do not yet have convergence of the derivatives, we can pass to the limit in the distributional sense. Since $d\varphi_j\to0$ in $L^2$ it follows that $d\Phi(u)=0$. Therefore $\Phi(u)\in\Ccal(g,o)$. We can now apply the elliptic estimate to $u_j-u$.
Both its Dirac derivative and its $L^2$ norm tend to zero,
so this gives us strong $H^1$ convergence of the spinors.
To obtain the same convergence for the associated three-forms,
we use the parallel bilinear map $B$ defining $\Phi$ and write
\[
 \nabla\Phi(u_j)-\nabla\Phi(u)
 =2B(u_j,\nabla(u_j-u))+2B(u_j-u,\nabla u).
\]
The first term tends to zero in $L^2$ because $u_j$ is pointwise unit and $u_j\to u$ in $H^1$. The second does so because $\nabla u$ is bounded and $u_j\to u$ in $L^2$. This proves the claimed strong $H^1$ convergence. Finally, a zero infimum in the last assertion would give a sequence to which this compactness result applies, contradicting emptiness of $\Ccal(g,o)$.
\end{proof}

The quadratic model also allows us to estimate how far a compatible form is from the space of closed forms when $\|d\varphi\|_{L^{2}}$ is small. The argument uses the finite-dimensional semialgebraic \L ojasiewicz inequality \cite{BCR}, after projecting a defining spinor onto the Dirac kernel.

\begin{theorem}\label{thm:stability}
Suppose $\Ccal(g,o)\ne\varnothing$. There exist constants $C,\epsilon>0$ and $\theta\in(0,1]$, depending on $(M,g,o)$, such that
\[
 \operatorname{dist}_{H^1}\bigl(\varphi,\Ccal(g,o)\bigr)
 \leq C\|d\varphi\|_{L^2}^{\theta}
\]
for every $\varphi\in\mathcal P(g,o)$ with $\|d\varphi\|_{L^2}<\epsilon$.
\end{theorem}
\begin{proof}
We first work in one spin sector with nonempty closed locus. Let $u$ be a unit defining spinor, let $h$ be its $L^2$ projection onto $\Hcal_\sigma$, and put $\eta=u-h$. The elliptic estimate on the orthogonal complement of the kernel and Lemma \ref{lem:dirac-bound} give
\[
 \|\eta\|_{H^1}\leq C\|D_\sigma u\|_{L^2}
 \leq C\|d\varphi\|_{L^2}.
\]
The $L^2$ norm of $h$ is bounded by that of $u$. All norms on the finite-dimensional smooth space $\Hcal_\sigma$ are equivalent, so $h$ and $\nabla h$ are uniformly bounded pointwise. Since $u$ is unit, polarization gives
\[
 \bigl\||h|^2-1\bigr\|_{L^2}\leq C\|\eta\|_{L^2}.
\]
As in the last part of the proof of Theorem \ref{thm:gap},
\[
 \nabla\Phi(u)-\nabla\Phi(h)
 =2B(u,\nabla\eta)+2B(\eta,\nabla h).
\]
It follows that $\|d\Phi(h)\|_{L^2}\leq C\|d\varphi\|_{L^2}$. Thus both components of $Q(h,h)$ have this bound. Orthogonality also gives $(h,h)_0=1-(\eta,\eta)_0$, so for small defect we may normalize $w=h/(h,h)_0^{1/2}$. Then
\[
 (w,w)_0=1,\qquad
 \|Q(w,w)\|_{\mathrm{av}}\leq C\|d\varphi\|_{L^2}.
\]
On this finite-dimensional unit sphere the polynomial map $w\mapsto Q(w,w)$ has compact zero set. The semialgebraic \L ojasiewicz inequality gives a zero $v$ with
\[
 \|w-v\|_0\leq C\|Q(w,w)\|_{\mathrm{av}}^\theta.
\]
Equivalence of norms on $\Hcal_\sigma$ and the estimate for $\eta$ imply $\|u-v\|_{H^1}\leq C\|d\varphi\|_{L^2}^\theta$. Zeros $v$ are pointwise unit and have uniformly bounded derivatives. The same bilinear estimate therefore proves the assertion for $\Phi(u)$ and $\Phi(v)$.

To finish the proof, we pass from a single spin structure to the
whole fibre. There are only finitely many spin structures on $M$.
For each one that is induced by a compatible form but admits no
compatible closed form, the argument of Theorem \ref{thm:gap},
restricted to that spin structure, gives a positive lower bound
for $\|d\varphi\|_{L^2}$. We choose $\epsilon$ smaller than all
these bounds. Every compatible form with
$\|d\varphi\|_{L^2}<\epsilon$ therefore belongs to a spin structure
whose closed locus is nonempty, where the estimate above applies.
Since there are finitely many such spin structures, we can take
the largest of their constants $C$, the smallest of their
exponents $\theta$, and decrease $\epsilon$ further so that all
the estimates hold and $\epsilon\leq1$. This gives the claimed
estimate on the whole fibre.
\end{proof}

\subsection{Scalar-flat metrics and parallel spinors}\label{sec:scalarflat}
The fixed-metric space also underlies the isometric flow of Dwivedi--Gianniotis--Karigiannis \cite{DGK}. Its critical points have divergence-free full torsion, a variational condition different from the Dirac harmonicity used here \cite{Nied}. Related formulations use spinorial energy \cite{AWW}, special-metric flows \cite{WW}, fixed-metric torsion energy \cite{Bagaglini}, and harmonic sections of geometric-structure bundles \cite{LSE}. Grigorian \cite{GrigSurvey} surveys these approaches, and Garrone \cite{Garrone} give non closed examples with divergence-free torsion on solvable groups.

We use the following differential-form conventions \cite{Karigiannis,LW}. The Hodge star is $*=*_{g,o}$, the formal adjoint of $d$ is $d_g^*$, and the Hodge Laplacian is $\Delta_g=dd_g^*+d_g^*d$. On three-forms in dimension seven, $d_g^*=-*d*$. For a closed $G_2$-structure, its torsion two-form $\tau_2$ is specified by
\[
 d(*\varphi)=\tau_2\wedge\varphi,\qquad
 \tau_2\in\Omega^2_{14}(M),\qquad
 d_g^*\varphi=\tau_2,\qquad \Delta_g\varphi=d\tau_22.
\]
Here $\Omega^2_{14}(M)$ denotes the smooth sections of $\Lambda^2_{14}T^*M$. Equivalently, $*(\beta\wedge\varphi)=-\beta$ for $\beta\in\Omega^2_{14}(M)$. The last two identities follow from this algebraic identity and $d\varphi=0$. Moreover, $(d_g^*)^2=0$ gives $d_g^*\tau_2=0$. In the full-torsion convention $T_{ij}=-\tfrac12(\tau_2)_{ij}$ for closed structures \cite{LW}, this is the divergence-free condition $\nabla^iT_{ij}=0$, with indices raised using $g$. hence every compatible closed form is a stationary point of the isometric flow \cite{DGK}. This is also a special case of \cite{GrigorianHeat}, which gives $\operatorname{div}T=0$ whenever $\tau_0$ is constant and $\tau_1=0$. 

There is also a useful integrated identity for every compatible form,
without assuming that it is closed. Let $\varphi=\Phi(u)\in\mathcal P(g,o)$,
where $u$ is its unit defining spinor in the induced spin structure
$\sigma$. The relation $\nabla_Xu=-T(X)\cdot u/2$ gives
$|\nabla u|^2=|T|^2/4$. Since $M$ is compact and has no boundary,
integrating the Schr\"odinger--Lichnerowicz formula and using
Lemma \ref{lem:dirac-bound} gives
\begin{equation}\label{eq:torsion-energy}
 \begin{aligned}
 \int_M |T|^2\,d\vol_g
 &=4\|D_\sigma u\|_{L^2(g)}^2-\int_M\Scal(g)\,d\vol_g\\
 &\leq\|d\varphi\|_{L^2(g)}^2-\int_M\Scal(g)\,d\vol_g.
 \end{aligned}
\end{equation}
The scalar-curvature integral is fixed by $g$, so this gives us a bound
for the full torsion energy whenever $\|d\varphi\|_{L^2(g)}$ is bounded.
This explains why the missing torsion component in the pointwise
closedness equation does not prevent an integrated energy bound.

The normalization of the defining spinor also gives a useful pointwise identity. If $|u|=1$ and $D_\sigma u=0$, the Schr\"odinger--Lichnerowicz formula implies
\begin{equation}
 |\nabla u|^2=-\tfrac14\Scal(g),\qquad |T|^2=-\Scal(g).
\end{equation}
Indeed, pairing $\nabla^*\nabla u=-\Scal(g)u/4$ with $u$ and differentiating $|u|^2=1$ gives the first equality pointwise. The second follows from $\nabla_Xu=-T(X)\cdot u/2$. For a closed structure, $T_{ij}=-(\tau_2)_{ij}/2$ and the tensor norm is $|T|^2=|\tau_2|^2/2$. This recovers the scalar formula below. It also gives a direct first-derivative bound for compatible closed forms whenever the metric and its scalar curvature are bounded. Higher regularity in Theorem \ref{thm:proper} is supplied by the Dirac equation.

 Bryant's scalar-curvature formula \cite[equation (4.36)]{Bryant} is
\begin{equation}\label{eq:scalar}
 \Scal(g)=-\tfrac12|\tau_2|_g^2.
\end{equation}
A positive scalar-curvature value anywhere therefore excludes a compatible closed form. In particular, $\mathscr A_o\subset\{g:\Scal(g)\leq0\}$. If $\Scal(g)\equiv0$, every such form is torsion-free. These statements require neither compactness nor completeness. All compatible closed forms have the same pointwise torsion norm. A related restriction is Cleyton--Ivanov's theorem \cite{CI}: on a compact manifold without boundary, an Einstein metric induced by a closed $G_2$-structure forces that structure to be torsion-free. Ball's classification with confromally flat induced metric \cite{BallFlat} gives a further precedent for imposing metric restrictions on closed structures. His quadratic torsion condition \cite{Ball} is different from the quadratic spinor equations used here. Xu--Zheng \cite{XZ} study the space of closed forms in a fixed de Rham class, allowing the metric to vary.

Suppose now that compact connected $M$ has a torsion free form $\varphi_0$ inducing $(g,o,\sigma_0)$, with parallel unit spinor $u_0$. The metric is Ricci-flat by Bonan's theorem \cite{Bonan}. The integrated Schr\"odinger--Lichnerowicz identity
\[
 D_{\sigma_0}^2=\nabla^*\nabla+\tfrac14\Scal(g)
\]
shows that harmonic spinors are exactly parallel spinors, compare \cite[Remark~9.7]{Grigorian}. Their pairwise inner products and bilinear three-forms are parallel, so $Q$ vanishes identically on $\Hcal_{\sigma_0}$. The parallel splitting
\[
 S_{\sigma_0}=\R u_0\oplus\{X\cdot u_0:X\in TM\}
\]
identifies parallel spinors with constants and parallel vector fields. By the integrated Bochner identity, parallel one forms on compact Ricci-flat $M$ are exactly harmonic one-forms. Hence
\begin{equation}\label{eq:lin}
 \dim_\R\Hcal_{\sigma_0}=1+b_1(M),\qquad
 \Ccal_{\sigma_0}(g,o)\cong\RP^{b_1(M)}.
\end{equation}
The parallel-vector-field argument is the one used by Grigorian
\cite[Theorem 7.6]{Grigorian} and Lin
\cite[Propositions 3.1--3.2]{Lin}. Their statements describe the whole
metric class, while their proofs use a globally defined octonion section
or a global pair $(f,\omega)$ relative to a reference form. Such a lift
parametrizes the spin sector of that reference form. A section of the
projective bundle need not admit this global lift, so the single
$\RP^{b_1(M)}$ description in those statements needs this sector
qualification. The following proposition records the full fibre.

\begin{proposition}\label{prop:parallel-sectors}
Let $(M,g,o)$ be a compact connected oriented seven-manifold without
boundary, with $\Scal(g)=0$. Write $\mathfrak S_{\parallel}(g,o)$ for
the set of spin structures carrying a nonzero parallel spinor. Then
\begin{equation}\label{eq:all-parallel}
 \Ccal(g,o)\cong
 \coprod_{\sigma\in\mathfrak S_{\parallel}(g,o)}\RP^{b_1(M)}.
\end{equation}
Every form in this fibre is torsion-free. The disjoint union is empty
when no spin structure carries a nonzero parallel spinor.
\end{proposition}
\begin{proof}
Any closed form in any spin sector is torsion-free by
\eqref{eq:scalar}. Conversely, a nonzero parallel spinor has constant
positive length and, after normalization, defines a torsion-free form.
In each nonempty sector, choose a parallel unit spinor and repeat the
splitting and Bochner argument above. It gives $Q\equiv0$ on
$\Hcal_\sigma$ and $\dim\Hcal_\sigma=1+b_1(M)$, so Theorem
\ref{thm:main} identifies that sector with $\RP^{b_1(M)}$.
Corollary \ref{cor:all} then gives the stated disjoint union.
\end{proof}

The distinction is visible in the compact flat example
\eqref{eq:flat-two}, where $b_1(M)=3$ and the full fibre is
$\RP^3\sqcup\RP^3$. Each component belongs to a different spin
structure. Thus the full fibre need not be the single projective
space obtained from a global lift relative to one reference form.

For a connected Euclidean domain $U\subset\R^7$, \eqref{eq:scalar} forces every compatible closed form to be parallel, hence constant in Cartesian coordinates. Conversely every constant compatible form is closed. Thus the complete fibre is $SO(7)/G_2\cong\RP^7$.

\section{The invariant fibre on almost-abelian Lie groups}\label{sec:invariant}
\subsection{The invariant closure equation}
Let $(V^6,\langle\ ,\ \rangle,o_V)$ be an oriented Euclidean space, and let
\[
 \mathfrak g_A=V\oplus\R\nu,\qquad
 [\nu,X]=AX,\quad [X,Y]=0\quad(X,Y\in V).
\]
The metric makes $\nu$ unit and orthogonal to $V$, and the orientation is $o_V\wedge\nu^*$. No compactness assumption is made on the associated simply connected Lie group. Put
\[
 \Ccal_A=\{\varphi\in\Lambda^3\mathfrak g_A^*:g_\varphi=g,\ o_\varphi=o,
                       \ d_{\mathfrak gA}\varphi=0\}.
\]
This is exactly the space of its left-invariant compatible closed forms.

For $\beta\in\Lambda^kV^*$ define
\[
 (\theta(A)\beta)(X_1,\ldots,X_k)
       =-\sum_{j=1}^k\beta(X_1,\ldots,AX_j,\ldots,X_k).
\]
With this notation, the Chevalley--Eilenberg differential is given by
\begin{equation}\label{eq:ce}
 d\beta=\nu^*\wedge\theta(A)\beta,\qquad d\nu^*=0.
\end{equation}
This is the formula in \cite[Proposition 2.17(b)]{Freibert},
and it also appears in \cite{FreibertCocalibrated,Lauret,Moreno}.
To check the sign with our convention, evaluate $d\beta$ on
$(\nu,X_1,\ldots,X_k)$. The brackets $[\nu,X_j]=AX_j$ give
exactly the sum defining $\theta(A)\beta$, while all brackets
among vectors in $V$ vanish.

Every positive form inducing this metric and orientation has a unique expression
\begin{equation}\label{eq:su3}
 \varphi=\omega_J\wedge\nu^*+\operatorname{Re}\Omega,
 \qquad \omega_J(X,Y)=\langle JX,Y\rangle,
\end{equation}
where $J$ is an orthogonal complex structure inducing $o_V$, and $\Omega$ is a complex $(3,0)$-form normalized by
\[
 \frac{i}{8}\Omega\wedge\overline{\Omega}=\vol_V
       =\frac16\omega_J^3.
\]
This follows from the transitive action of $G_2$ on unit vectors, with stabilizer $SU(3)$, so the space of such pairs is $SO(6)/SU(3)$.

\begin{lemma}\label{lem:closure}
The form \eqref{eq:su3} is closed if and only if
\begin{equation}\label{eq:commute}
 AJ=JA,\qquad \tr_{\C,J}A=0.
\end{equation}
Write $o_J$ for the complex orientation of $(V,J)$. Then $\Ccal_A$ is the circle bundle of normalized complex volume forms over
\[
 \mathcal J_A=\{J:J^2=-\Id,\ J^tJ=\Id,\ o_J=o_V,\ [A,J]=0,
                         \tr_{\C,J}A=0\}.
\]
\end{lemma}
\begin{proof}
This is the closure calcuilation in Freibert's criterion \cite[Theorem~3.1(a)]{Freibert}, written with the metric fixed. Equation \eqref{eq:ce} gives $d\varphi=\nu^*\wedge\theta(A)\operatorname{Re}\Omega$. The Lie algebra stabilizer of $\operatorname{Re}\Omega$ is $\mathfrak{sl}(3,\C)$ by Hitchin's stable-form construction \cite{Hitchin}, see also \cite[Lemma~2.7(i)]{Freibert}. Explicitly, $\operatorname{Re}\Omega$ determines the pair $\{J,-J\}$. If $\theta(A)\operatorname{Re}\Omega=0$, the one-parameter group $\exp(tA)$ preserves this pair. Continuity in $t$, together with $\exp(0A)=\Id$, forces it to preserve $J$ itself. Differentiating gives $[A,J]=0$. For a complex-linear endomorphism,
\[
 \theta(A)\Omega=-(\tr_{\C,J}A)\Omega.
\]
The real and imaginary parts of $\Omega$ are linearly independent, so the real part of the right-hand side vanishes exactly when the complex trace is zero. The converse follows from the same identity. For each admissible $J$, normalized volume forms form a circle, and \eqref{eq:su3} reconstructs the form uniquely.
\end{proof}

The classification below retains the individual forms, including those related by orthogonal automorphisms.

\subsection{The orthogonal centralizer without normality}
The closure criterion does not require $A$ to be normal. Before using spectral theory, we can therefore describe the connected components for an arbitrary endomorphism in terms of its orthogonal centralizer. Let
\[
 \mathscr K_A=\{R\in SO(V):RA=AR\}^{0},
\]
where the superscript denotes the identity component. For $J\in\mathcal J_A$ define
\[
 \mathscr L_J=\{R\in\mathscr K_A:RJ=JR\},\qquad
 \mathscr H_J=\ker\bigl(\det_{\C,J}:\mathscr L_J\longrightarrow S^1\bigr).
\]
These definitions distinguish the compact symmetry group from the simply connected almost-abelian Lie group itself. The action on the geometric data is
\[
 R\cdot(J,\Omega)=\bigl(RJR^{-1},(R^{-1})^*\Omega\bigr).
\]
It preserves the prescribed metric, orientation and closure equation.

\begin{theorem}\label{thm:centralizer}
For every endomorphism $A$ of an oriented Euclidean six-space, the space $\Ccal_A$ is either empty or a disjoint union of at most four connected compact homogeneous spaces. Each component is of the form $\mathscr K_A/\mathscr H_J$, for an admissible $J$, and the projection to its complex structure is the principal circle bundle
\[
 \mathscr K_A/\mathscr H_J\longrightarrow\mathscr K_A/\mathscr L_J.
\]
In particular, every nonempty component has positive dimension. The bound of four components is attained.
\end{theorem}
\begin{proof}
If $J$ is an orthogonal complex structure commuting with $A$, it is skew-adjoint and belongs to the Lie algebra of $\mathscr K_A$. The conjugacy theorem for maximal tori in connected compact Lie groups \cite{BtD} puts it, by conjugation in $\mathscr K_A$, in the Lie algebra $\mathfrak t$ of a fixed maximal torus of $\mathscr K_A$. Embed this torus in a maximal torus of $SO(6)$. In the associated orthogonal decomposition of $V$ into three oriented two-planes, an element of it Lie algebra is
\[
 \operatorname{diag}(a_1J_0,a_2J_0,a_3J_0),
 \qquad J_0=\begin{pmatrix}0&-1\\1&0\end{pmatrix}.
\]
Its square is $-\Id$ exactly when each $a_i$ is $1$ or $-1$. Exactly four of these eight choices have the prescribed complex orientation. Some need not lie in $\mathfrak t$, and some can be conjugates under $\mathscr K_A$, but every admissible $J$ is conjugate to one of the remaining choices. The condition $\tr_{\C,J}A=0$ is constant on a conjugacy orbit. Hence $\mathcal J_A$ is a union of at most four compact orbits.

For each such $J$, the group $\mathscr L_J$ acts on normalized complex volume forms through the inverse complex determinant. This character is surjective: $\exp(tJ)$ lies in $\mathscr L_J$, and
\[
 \det_{\C,J}\exp(tJ)=e^{3it}.
\]
Hence $\mathscr K_A$ is transitive on all pairs above this orbit, with stabilizer $\mathscr H_J$, and the quotient of $\mathscr L_J$ by this stabilizer is a circle. Since $\mathscr K_A$ is connected, the total space is connected. Distinct such orbits are disjoint compact subsets, and there are finitely many of them, so they are precisely the connected components. The symmetric $2+2+2$ case of Corrollary \ref{cor:symmetric} attains four components.
\end{proof}

The homogeneous-space description shows that every component of
the invariant fibre is smooth. For the general compact fibre,
the algebraic description in Section \ref{sec:compact} does not
by itself establish smoothness. When $A$ is normal, we can
describe its centralizer and the corresponding homogeneous spaces
explicitly using its spectral decomposition. Before doing this,
we consider the non-normal case. In real dimension six, the
failure of normality severely restricts the orthogonal complex
structures that commute with $A$, as the following lemma shows.

\begin{lemma}\label{lem:two-J}
Suppose an endomorphism $A$ of a Euclidean six-space comutes with two orthogonal complex structures $J$ and $J'$ which do not commute with each other. Then $A$ is normal.
\end{lemma}
\begin{proof}
The self-adjoint endomorphism $P=-\frac12(JJ'+J'J)$ commutes with $A,J,J'$. For a unit vector $x$, $\langle Px,x\rangle=\langle Jx,J'x\rangle$, so its eigenvalues belong to $[-1,1]$. On its $c$-eigenspace with $|c|<1$, put
\[
 I=\frac{J'-cJ}{\sqrt{1-c^2}}.
\]
The relation $JJ'+J'J=-2c\Id$ on that eigenspace gives $I^2=-\Id$, $I^t=-I$ and $JI=-IJ$. Thus $J,I$ define an orthogonal quaternionic structure there, and the dimension of this eigenspace is a multiple of four. On the $c=1$ and $c=-1$ eigenspaces one has respectively $J'=J$ and $J'=-J$: this follows by computing the squared norm of $(J'\mp J)x$. If there were no eigenvalue in $(-1,1)$, $J$ and $J'$ would consequently commute. There is therefore such an eigenspace, and in dimension six it has dimension four, with an orthogonal two-dimensional complement.

Both summands are preserved by $A$. On the four-dimensional summand, $A$ commutes with $J$ and $I$ and is therefore a quaternionic-linear endomorphism of a quaternionic line. Upon identifying this line isometrically with $\mathbb H$, it is multiplication on the opposite side by a quaternion $q$. Its adjoint is multiplication by $\overline q$, so it commutes with its adjoint. On the two-dimensional complement, $A$ commutes with the orthogonal complex structure $J$ and hence has the form $a\Id+bJ$, which is also normal. The orthogonal direct sum is normal.
\end{proof}

\begin{theorem}\label{thm:nonnormal}
If $A$ is not normal, then $\Ccal_A$ is empty, a circle, or a disjoint union of two circles. Both nonempty possibilities occur for nilpotent $A$.
\end{theorem}
\begin{proof}
By Lemma \ref{lem:two-J}, any two orthogonal complex structures commuting with $A$ commute with each other. Suppose three distinct such structures $J_1,J_2,J_3$ induce the same orientation. Viewing $J_2$ and $J_3$ as commuting skew-Hermitian endomorphisms of $(V,J_1)$, simultaneous unitary diagonalization gives an orthogonal decomposition into three real two-planes on which each $J_i$ is $J_0$ or $-J_0$. Orient each plane by $J_1$. Its sign triple is then $(1,1,1)$, and the other triples have product $+1$ because their complex orientations agree. The four possibilities are
\[
 (1,1,1),\quad(1,-1,-1),\quad(-1,1,-1),\quad(-1,-1,1),
\]
and any three of these triples have determinant $4$ or $-4$. It follows that the three operators which equal $J_0$ on one two-plane and zero on the others are linear combinations of $J_1,J_2,J_3$. They commute with $A$, so their squares show that $A$ preserves each two-plane. On each plane it is complex scalar, and hence normal. This contradicts the hypothesis. There are consequently at most two orthogonal complex structures of the prescribed orientation commuting with $A$, before imposing the trace condition. Lemma \ref{lem:closure} supplies one circle over each admissible one.

For completeness, let $N_r$ be the real nilpotent Jordan block of size $r$, with ones directly above the diagonal, and use the standard tensor-product metrics. For
\[
 A=N_3\otimes I_2\quad\hbox{on }\R^3\otimes\R^2,
\]
the matrices $N_3,N_3^t$ generate the full algebra $M_3(\R)$. Their common commutant is therefore $I_3\otimes M_2(\R)$. An orthogonal complex structure commuting with $A$ also commutes with $A^t$, by transposing the commutation relation. It must thus be $I_3\otimes J_0$ or $-I_3\otimes J_0$. These have opposite orientations, so exactly one is allowed. Its complex trace is zero because $A$ is nilpotent.. This gives a single circle.

For
\[
 A=(N_2\otimes I_2)\oplus0_{\R^2},
\]
the common kernel of $A$ and $A^t$ is the last two-plane, which every commuting orthogonal complex structure preserves. On the four-dimensional summand the common commutant is $I_2\otimes M_2(\R)$, giving the two choices $\pm(I_2\otimes J_0)$, both of the same orientation. On the remaining two-plane the orientation condition fixes one of $\pm J_0$. Thus there are exactly two admissible complex structures, again with trace zero, giving two circles. Both displayed matrices are non-normal since $N_rN_r^t\ne N_r^tN_r$ for $r\geq2$.
\end{proof}

More generally, the theorem shows that a positive-dimensional family of admissible complex structures in six dimensions forces normality of $A$. The circle of complex volume forms persisits without normality, but the additional freedom in the complex structure does not.

To distinguish the non-normal alternatives, we first compute the skew-adjoint endomorphisms commuting with $A$. This gives us a linear space in which the admissible complex structures are determined by a finite system of quadratic equations. Let
\[
 \mathfrak c_A=\{B\in\mathfrak{so}(V):BA=AB\},
\]
and choose a real basis $B_1,\ldots,B_r$. Transposing the commutation relation shows that every member also commutes with $A^t$, so this is the skew-adjoint part of the commutant of the real $*$-algebra generated by $A$. Define
\begin{equation}\label{eq:commutant-test}
 J(x)=\sum_{\alpha=1}^r x_\alpha B_\alpha,
 \qquad J(x)^2=-\Id,\quad\tr A=0,\quad\tr(J(x)A)=0.
\end{equation}
For a commuting complex structure,
$\tr_{\C,J}A=\tfrac12\tr A-\tfrac{i}{2}\tr(JA)$.
Therefore, after retaining the solutions of the prescribed orientation, \eqref{eq:commutant-test} has exactly as many solutions as $\Ccal_A$ has circles when $A$ is non-normal. Theorem \ref{thm:nonnormal} bounds this number by two. In particular, one circle means that this explicit finite system has exactly one oriented solution.

For non-normal $A$, there is also a geometric description of when
the two-circle case occurs. It occurs precisely when there is an
orthogonal splitting $V=V_4\oplus V_2$ such that
\[
 A=A_4\oplus c\Id_{V_2},\qquad c\in\R,
 \qquad [A_4,J_4]=0,\qquad
 \tr_{\C,J_4}A_4=-c
\]
for some orthogonal complex structure $J_4$ on $V_4$.

To see this, suppose that $J$ and $J'$ are two distinct admissible
complex structures. They commute by Lemma \ref{lem:two-J}, so
$V$ splits orthogonally into the subspaces on which they agree
and disagree. Both subspaces are preserved by $A$. Since $J$
and $J'$ induce the same orientation, the complex dimension of
the disagreement subspace is even. They are distinct and
$\dim_{\C}V=3$, so the disagreement subspace has real dimension
four. We denote it by $V_4$, and denote the agreement subspace
by $V_2$. Thus
\[
 J=J_4\oplus J_2,\qquad J'=-J_4\oplus J_2,
 \qquad A=A_4\oplus A_2.
\]
Put
\[
 z=\tr_{\C,J_4}A_4,\qquad w=\tr_{\C,J_2}A_2.
\]
The two complex trace conditions are
\[
 z+w=0,\qquad \overline z+w=0.
\]
Hence $z$ and $w$ are real. Since $A_2$ is a complex-linear
endomorphism of a complex line, its real complex trace implies
$A_2=c\Id_{V_2}$ for some $c\in\R$, and the first trace equation
gives $z=-c$.

Conversely, suppose that the displayed splitting exists. Choose
$J_2$ so that $J_4\oplus J_2$ induces the prescribed orientation.
Then $-J_4\oplus J_2$ induces the same orientation, since the
sign has been reversed on a complex two-plane. Both structures
commute with $A$, and both have zero complex trace. Lemma
\ref{lem:closure} therefore gives a circle of closed forms over
each of them. Theorem \ref{thm:nonnormal} shows that these are
the only two circles.

\subsection{Spectral decomposition and homogeneous components}
Assume $AA^t=A^tA$. Its orthogonal real spectral decomposition is
\begin{equation}\label{eq:spectral}
 V=\bigoplus_{\lambda\in\mathcal R}V_\lambda
       \ \oplus\!\bigoplus_{z=a_z+ib_z\in\mathcal Z}W_z,
 \qquad b_z>0,
\end{equation}
where $\mathcal R=\operatorname{spec}(A)\cap\R$, $A|_{V_\lambda}=\lambda\Id$ and
$A|_{W_z}=a_z\Id+b_zI_z$. Here $I_z^2=-\Id$, $I_z$ is orthogonal, and $\dim_\C(W_z,I_z)=m_z$. Distinct elements of $\mathcal Z$ are distinct complex eigenvalues in the upper half-plane. If any $\dim_\R V_\lambda$ is odd, set the admissible-index set below to be empty. Otherwise write $\dim_\R V_\lambda=2m_\lambda$, so that
\[
 \sum_\lambda m_\lambda+\sum_zm_z=3.
\]
Choose a reference orientation on each real eigenspace. The complex structures $I_z$ orient the non-real blocks. Their product gives a reference orientation on $V$, whose comparison with $o_V$ is denoted by $\delta\in\{1,-1\}$: $o_V=\delta\,o_{\rm ref}$.

An index consists of signs $\epsilon_\lambda\in\{1,-1\}$ and integers $0\le p_z\le m_z$. Set $q_z=m_z-p_z$. Call this index admissible when
\begin{align}
 \sum_\lambda m_\lambda\lambda+\sum_zm_za_z&=0,\label{eq:realtrace}\\
 \sum_zb_z(p_z-q_z)&=0,\label{eq:imagtrace}\\
 \prod_\lambda\epsilon_\lambda\,(-1)^{\sum_zq_z}&=\delta.\label{eq:orientation}
\end{align}
An empty product is one. Equation \eqref{eq:realtrace} is $\tr_\R A=0$ divided by two, the unimodularity condition recorded in \cite[Remark~3.3]{Freibert}.

For one admissible index choose a reference $J$ whose complex orientation on $V_\lambda$ has sign $\epsilon_\lambda$, and which is $I_z$ on a complex $p_z$-plane and $-I_z$ on its orthogonal complement. For normal $A$, the connected centralizer introduced above is
\begin{equation}
 \mathscr K_A=\prod_\lambda SO(2m_\lambda)\times\prod_zU(m_z).
\end{equation}
Its stabilizer of $J$ is the subgroup
\[
 \mathscr L_J=\prod_\lambda U(m_\lambda)\times
                   \prod_z\bigl(U(p_z)\times U(q_z)\bigr).
\]
The embeddings on the real blocks use the chosen $J$ and its orientation component. On non-real blocks the unitary groups use $I_z$. Omit factors $U(0)$. The complex determinant for the total complex structure $J$ is the character
\begin{equation}
 \chi_J((R_\lambda),(P_z,Q_z))
   =\prod_\lambda\det_\C R_\lambda
            \prod_z\det_\C P_z\,(\det_\C Q_z)^{-1},
 \qquad \mathscr H_J=\ker\chi_J.
\end{equation}
The inverse determinant on a negative block comes from replacing $I_z$ by $-I_z$.

\begin{theorem}\label{thm:normal}
For a normal endomorphism $A$ of an oriented Euclidean six-space, the following hold.
\begin{enumerate}
\item The space $\Ccal_A$ is nonempty if and only if every real eigenspace has even dimension and an admissible index exists.
\item For each admissible index, there is exactly one connected component, and it is diffeomorphic to $\mathscr K_A/\mathscr H_J$. In particular,
\begin{equation}
 \Ccal_A\cong\coprod_{\text{admissible indices}}\mathscr K_A/\mathscr H_J.
\end{equation}
\item The component is a principal circle bundle over
\begin{equation}\label{eq:base}
 \prod_\lambda SO(2m_\lambda)/U(m_\lambda)
     \times\prod_z\operatorname{Gr}_{p_z}(\C^{m_z}),
\end{equation}
and has dimension
\begin{equation}\label{eq:dimension}
 1+\sum_\lambda m_\lambda(m_\lambda-1)+2\sum_zp_zq_z.
\end{equation}
\end{enumerate}
Every component is an orbit of the identity component of the orthogonal automorphisms fixing $\nu$ and commuting with $A$.
\end{theorem}
\begin{proof}
If $J$ commutes with $A$, it also commutes with $A^t$. Indeed, transpose $AJ=JA$ and use $J^t=-J$. Hence it preserves the joint spectral decomposition of $S=(A+A^t)/2$ and $A_-=(A--A^t)/2$, and therefore all blocks in \eqref{eq:spectral}.

On a real eigenspace, $J$ is an arbitrary orthogonal complex structure. Such a structure exists precisely in even dimension. For each of its two possible complex orientations the space is $SO(2m_\lambda)/U(m_\lambda)$. These are connected spaces. Their labels are the signs $\epsilon_\lambda$.

On $W_z$, commutation with $A$ is commutation with $I_z$. Thus $J$ is a complex-linear skew-Hermitian operator on $(W_z,I_z)$, with $J^2=-\Id$. The unitary spectral theorem decomposes $W_z$ orthogonally into its $+i$ and $-i$ eigenspaces of dimensions $p_z$ and $q_z$. Equivalently, $J=I_z$ on the first and $J=-I_z$ on the second. The choices form the connected complex Grassmannian $\operatorname{Gr}_{p_z}(\C^{m_z})$.

The complex trace of $A$ relative to this $J$ is
\begin{equation}
 \tr_{\C,J} A=
 \sum_\lambda m_\lambda\lambda+\sum_zm_za_z
           +i\sum_zb_z(p_z-q_z).
\end{equation}
On a negative non-real block the complex eigenvalue is conjugated, explaining the last sign. Reversing the complex structure on a complex $q_z$-plane changes its real orientation by $(-1)^{q_z}$. The total orientation is therefore exactly \eqref{eq:orientation}. Lemma \ref{lem:closure} now proves the existence assertion and the base description \eqref{eq:base}.

We now apply Theorem \ref{thm:centralizer} with the centralizer
and stabilizer computed above. It identifies the total space
over each admissible choice of indices with
$\mathscr K_A/\mathscr H_J$, and gives the circle bundle over
$\mathscr K_A/\mathscr L_J$. Each such total space is connected.
The indices $p_z$ and $\epsilon_\lambda$ are locally constant:
the former are dimensions of eigenspaces, and the latter are
orientation signs. Distinct admissible choices therefore give
distinct connected components. To compute their dimensions,
we use
\[
 \dim SO(2m)-\dim U(m)=m(m-1),\qquad
 \dim U(m)-\dim U(p)-\dim U(q)=2pq,
\]
and passage from $\mathscr L_J$ to the codimension-one subgroup $\mathscr H_J$ adds one dimension. This proves \eqref{eq:dimension}.

An element $R\in \mathscr K_A$ extends by $R\nu=\nu$ to an orthogonal Lie algebra automorphism, since $RA=AR$. These extensions give the final assertion.
\end{proof}

\begin{remark}
The orientation condition only selects half of the indices satisfying the trace conditions. Replacing $J$ by $-J$ conjugates $\tr_{\C,J}A$ and reverses the orientation of the complex three-space. In spectral notation this sends
\[
 p_z\longleftrightarrow q_z,\qquad
 \epsilon_\lambda\longmapsto(-1)^{m_\lambda}\epsilon_\lambda.
\]
The left-hand side of \eqref{eq:orientation} changes by
$(-1)^{\sum_\lambda m_\lambda+\sum_zm_z}=-1$, while both trace equations remain satisfied. Thus $\Ccal_A$ is nonempty exactly when all real eigenspaces are even-dimensional and \eqref{eq:realtrace}--\eqref{eq:imagtrace} are solvable. More directly,
\[
 \Ccal_A(g,o)\longrightarrow\Ccal_A(g,-o),\qquad
 \varphi\longmapsto-\varphi
\]
is a diffeomorphism, so the two orientations have identical component counts and topology. On the $SU(3)$ data it is $(J,\Omega)\mapsto(-J,\Omega')$, where $\Omega'=-\overline{\Omega}$ is minus the complex conjugate of $\Omega$. Complex conjugation changes type $(3,0)$ for $J$ to type $(3,0)$ for $-J$, and
\[
 \operatorname{Re}\Omega'=-\operatorname{Re}\Omega,\qquad
 \frac{i}{8}\Omega'\wedge\overline{\Omega'}
 =\frac{i}{8}\overline{\Omega}\wedge\Omega=-\vol_V.
\]
\end{remark}

\begin{corollary}\label{cor:symmetric}
Let $A=A^t$. A compatible invariant closed form exists precisely when $\tr A=0$ and every eigenspace of $A$ is even-dimensional. In that case the complete fibre is
\[
\begin{array}{@{}lc@{}}\toprule
 \text{Spectral multiplicities of }A&\Ccal_A\\ \midrule
 2+2+2\ \text{with three distinct eigenvalues}&\displaystyle\bigsqcup_{j=1}^{4}S^1\\
 4+2\ \text{with two distinct eigenvalues}&\RP^3\sqcup\RP^3\\
 6\ \text{(necessarily }A=0\text{)}&\RP^7\\ \bottomrule
\end{array}
\]
\end{corollary}
\begin{proof}
There are no non-real blocks, so the imaginary trace equation is vacuous. For $r$ real eigenspaces, exactly $2^{r-1}$ choices of orientations satisfy \eqref{eq:orientation}. In the $2+2+2$ case each base factor is a point, so each of the four components is a circle.

In the $4+2$ case there are two components, each
\[
 (SO(4)\times SO(2))/S(U(2)\times U(1)).
\]
The determinant condition allows elimination of the $SO(2)$ factor, identifying this quotient with $SO(4)/SU(2)$. Under the double covering
$\operatorname{Spin}(4)=SU(2)\times SU(2)\to SO(4)$, the relevant $SU(2)$ stabilizer is one factor. The residual quotient is $SU(2)/\{\pm1\}=SO(3)\cong\RP^3$. The opposite complex orientation exchanges the two spin factors and gives the same quotient.

In the remaining case $A=\lambda\Id$ and trace zero forces $\lambda=0$. The fibre is $SO(6)/SU(3)$. Since $\operatorname{Spin}(6)=SU(4)$ and $SU(4)/SU(3)=S^7$, passage to $SO(6)$ identifies antipodal points, yielding $\RP^7$.
\end{proof}

\begin{corollary}\label{cor:table}
Every nonempty component of $\Ccal_A$ is a projective linear subspace of the space $\RP^7$ of invariant compatible forms. Its dimension is one, three or seven. The number of components is $0$, $1$, $2$ or $4$.

For normal $A$ the complete list is Table \ref{tab:normal}. In that table, $z,w,z_i$ lie in the upper half-plane, all indicated real eigenvalues are real, and multiplicities are multiplicities in the complexified real six-dimensional representation.
\end{corollary}

\begin{table}[htbp]
\centering\small
\renewcommand{\arraystretch}{1.5}
\begin{tabular}{@{}p{.36\textwidth}p{.39\textwidth}l@{}}
\toprule
Spectrum & Condition & $\Ccal_A$\\\midrule
$0$ (multiplicity $6$) & & $\RP^7$\\
$\lambda$ ($4$), $\mu$ ($2$) & $\mu=-2\lambda\ne0$ & $\RP^3\sqcup\RP^3$\\
$\lambda_1,\lambda_2,\lambda_3$ ($2$ each), distinct & $\lambda_1+\lambda_2+\lambda_3=0$ & $\coprod_{j=1}^4\RP^1$\\
$z,\overline z$ ($2$ each), $\mu$ ($2$) & $\mu=-2\operatorname{Re}z$ & $\RP^3$\\
$z,\overline z$ ($2$ each), $w,\overline w$ ($1$ each) & $w=-2\overline z$ & $\RP^1$\\
$\mu$ ($2$), $z,\overline z,w,\overline w$ ($1$ each), $z\ne w$ & $\operatorname{Im}z=\operatorname{Im}w$, $\mu+\operatorname{Re}z+\operatorname{Re}w=0$ & $\RP^1\sqcup\RP^1$\\
$z_1,z_2,z_3$ and conjugates ($1$ each), distinct & $z_1+z_2+\overline z_3=0$ after relabelling & $\RP^1$\\
All other spectra & & $\varnothing$\\\bottomrule
\end{tabular}
\caption{The normal invariant fibres.}\label{tab:normal}
\end{table}

\begin{proof}
We first obtain the list by applying the trace and orientation
conditions in Theorem \ref{thm:normal}. Since the total complex
block multiplicity is three, there are only a few possible
partitions to consider. A single block of multiplicity three must be real and zero, since $b(2p-3)=0$ has no solution for a non-real block. The other partitions are $2+1$ and $1+1+1$. For $2+1$, a non-real double block contributes imaginary trace $0$ or $\pm2\operatorname{Im}z$. A real remaining block requires the first choice, giving row four. A non-real remaining block requires the second, giving row five. A real double block and non-real single block cannot satisfy the imaginary trace equation. For $1+1+1$, three real blocks give row three. One real and two non-real blocks require opposite imaginary signs of equal magnitude, giving row six. Three non-real blocks require one imaginary part to be the sum of the other two, giving row seven. The real trace and orientation equations supply exactly the conditions and component counts in the table.

To identify the embeddings as projective linear spaces, begin with symmetric $S$. In the $2+2+2$ case each admissible $J$ is fixed. Its circle of volume forms is the projectivization of the real two-plane fixed by $SU(3)$ in the spin representation. In the $4+2$ case let $E$ be the two-dimensional eigenspace and choose an oriented orthonormal basis $e,f$ of $E$. The two components are characterized by
\[
 \varphi(\nu,e,f)=1\quad\hbox{or}\quad\varphi(\nu,e,f)=-1.
\]
Indeed, either equality says that $J$ preserves the oriented plane $E$, with the specified sign. On $E^\perp\cap V$ its compatible complex structure is then arbitrary of the required orientation. The symmetric Clifford involution $\nu\cdot e\cdot f$ has eigenspaces of real dimension four. Since
$\varphi(\nu,e,f)=-\langle\nu\cdot e\cdot f\,u,u\rangle$,
the two equalities describe their projectivizations. For $S=0$ the full spinor space gives $\RP^7$.

Now let $A=S+A_-$ be normal. The closure condition is the condition for $S$ together with preservation of the form by the one-parameter orthogonal group $\exp(tA_-)$. This follows by separating the real and imaginary parts of the complex trace, as in Subsection \ref{sec:torsion}. Because $[S,A_-]=0$, this group preserves each of the finitely many projective linear components $\mathbb P(W)$ for $S$. Its spin lift therefore preserves $W$. A line in $W$ is fixed precisely when its vector is annihilated by the skew-adjoint infinitesimal spin action $\widehat A_-$: a real invariant line of a skew-adjoint map has eigenvalue zero. Hence the closed part of that component is
\[
 \mathbb P\bigl(W\cap\ker\widehat A_-\bigr).
\]
The dimensions in Theorem \ref{thm:normal} now identify the three-dimensional entries as $\RP^3$, and the one-dimensional entries as $\RP^1$. For non-normal $A$, each component has a fixed $J$ by Theorem \ref{thm:nonnormal}, so the same $SU(3)$-fixed-plane argument applies. This proves the linearity assertion in all cases and completes the count.
\end{proof}

There is also a finite answer after quotienting by orthogonal Lie algebra automorphisms. Here only orientation-preserving automorphisms act on the prescribed fibre. For $A\ne0$ with nonempty fibre, the abelian ideal $V$ is unique: a second abelian hyperplane would meet $V$ in dimension five and would force $A$ to vanish on that intersection, whereas a nonzero complex-linear $A$ has real rank at least two. Therefore every orthogonal automorphism has the form
\begin{equation}\label{eq:aut}
 F|_V=R,\quad F\nu=\varepsilon\nu,
 \qquad RA=\varepsilon AR,\quad\det R=\varepsilon,\quad\varepsilon\in\{1,-1\}.
\end{equation}
It sends $J$ to $\varepsilon RJR^{-1}$. By Theorem \ref{thm:centralizer} the connected centralizer already acts transitively on each component. The quotient is therefore computed by the permutation of components induced by \eqref{eq:aut}.

For normal $A$ every nonempty row of Table \ref{tab:normal} gives one equivalence class, except row six, which gives two when $\mu\ne0$ and one when $\mu=0$. To check this, only the rows with more than one component require attention. For symmetric $A$, reflections in two real eigenplanes preserve total orientation and exchange the sign choices, transitively. In row six, automorphisms with $\varepsilon=1$ preserve each non-real spectral block and cannot exchange the two imaginary sign choices. An automorphism with $\varepsilon=-1$ requires $A$ and $-A$ to have the same spectrum. The unique real eigenvalue forces $\mu=0$. Conversely, then $w=-\overline z$, and conjugate-linear isometries exchanging the two non-real planes, together with a reflection on the real plane, give an $R$ of determinant $-1$ satisfying \eqref{eq:aut}. Its induced action exchanges the two components. For general background on automorphisms of closed structures, see Podest\`a--Raffero \cite{PR}.

Finally, the same closure calculation describes which metrics on a fixed almost-abelian Lie algebra are admissible. Fix $[f_7,X]=A_0X$ on $V$. Every metric is uniquely specified by a positive definite inner product $G$ on $V$, a vector $v\in V$, and $\ell>0$, by declaring $\nu=(f_7-v)/\ell$ unit and perpendicular to $V$. Its defining endomorphism is $A_0/\ell$. It admits a compatible closed form precisely when there is a complex structure $J$ with
\[
 [A_0,J]=0,\qquad\tr_{\C,J}A_0=0,\qquad
 G(JX,JY)=G(X,Y).
\]
Orientation can be corrected by replacing $J$ with $-J$. Thus the admissible inner products on $V$ form the union, over these complex structures, of their cones of positive definite Hermitian inner products. This follows directly from Lemma \ref{lem:closure}. It makes explicit why an algebraic existence criterion allowing arbitrary changes of basis need not imply existence for every prescribed metric.

\subsection{Torsion, harmonic spinors and compact quotients}\label{sec:torsion}
Write $A=S+A_-$, where $S^t=S$ and $A_-^t=-A_-$. Normality is equivalent to $[S,A_-]=0$. For an admissible $J$, both $S$ and $A_-$ commute with $J$. The trace condition in \eqref{eq:commute} separates into
\[
 \tr_{\C,J}S=0,\qquad \tr_{\C,J}A_-=0,
\]
since the first trace is real and the second purely imaginary. Hence, $A_-\in\mathfrak{su}(V,J)$ and $S$ is a trace-free Hermitian endomorphism. This is Moreno's closed-class splitting $\mathfrak{sym}^0_+\oplus\mathfrak{su}(3)$ \cite{Moreno}, and normality adds $[S,A_-]=0$. For any symmetric endomorphism $B$ commuting with $J$, write
\[
 \omega_B(X,Y)=\langle JBX,Y\rangle.
\]
In particular, $\omega_S$ is a primitive real $(1,1)$form.

The identities below follow from the closed almost-abelian
formulas of Lauret \cite{Lauret}. The torsion identities can
also be obtained by specializing Moreno's general formulas
\cite{Moreno} to the closed class, and the Ricci formula is the
unimodular case of Arroyo's almost-abelian formula \cite{Arroyo}.
We include the calculation to make the signs explicit in our
conventions and to keep track of the commutator terms when $A$
is not normal. For curvature formulas expressed in terms of
the intrinsic torsion of a six-dimensional $SU(3)$-structure,
see Bedulli--Vezzoni \cite{BV}. Here we compute the
seven-dimensional Ricci tensor directly from the
Levi--Civita connection of the almost-abelian metric.

\begin{proposition}\label{prop:torsion}
For every $\varphi\in\Ccal_A$, its torsion two-form, normalized by $d(*\varphi)=\tau_2\wedge\varphi$, satisfies
\begin{equation}\label{eq:torsion}
 \tau_2=-2\omega_S,\qquad
 d\tau_2=\nu^*\wedge(4\omega_{S^2}+2\omega_{[S,A_-]}),\qquad
 |\tau_2|^2=2\tr(S^2).
\end{equation}
The fixed metric has
\begin{equation}\label{eq:ricci}
 \operatorname{Ric}|_{V\times V}=\tfrac12[A,A^t]=-[S,A_-],\quad
 \operatorname{Ric}(\nu,V)=0,\quad
 \operatorname{Ric}(\nu,\nu)=-\tr(S^2).
\end{equation}
Here endomorphisms are identified with symmetric bilinear forms using the metric. In particular, $A$ is normal exactly when $\operatorname{Ric}|_{V	imes V}=0$. All these forms are torsion-free exactly when $S=0$. Otherwise every one is strictly closed and the scalar curvature is negative.
\end{proposition}
\begin{proof}
Set $\rho_- =\operatorname{Im}\Omega$. The Hodge dual in \eqref{eq:su3} is
\[
 *\varphi=\tfrac12\omega_J^2+\rho_-\wedge\nu^*.
\]
Because $A_-$ preserves $\omega_J$ and $S$ commutes with $J$,
\[
 \theta(A)\omega_J=-2\omega_S,
 \qquad d(*\varphi)=-2\nu^*\wedge\omega_S\wedge\omega_J.
\]
The primitive $(1,1)$ identity $*_6(\omega_S\wedge\omega_J)=-\omega_S$ therefore gives
\[
 d_g^*\varphi=-*_7d(*\varphi)=-2\omega_S.
\]
For a closed $G_2$-form, $d_g^*\varphi=\tau_2$, proving the first formula. Directly from the definition of $\theta$, one has $\theta(A)\omega_S=-2\omega_{S^2}-\omega_{[S,A_-]}$. Equation \eqref{eq:ce} proves the formula for $d\tau_2$. In a unitary eigenbasis, $S$ has real eigenvalues $s_1,s_2,s_3$, each repeated twice, and $\omega_S=\sum_js_j e^{2j-1,2j}$. Hence $|\omega_S|^2=\tfrac12\tr(S^2)$, giving the norm formula.

To compute the curvature, we begin with the Levi--Civita connection on invariant vector fields. The standard formula \cite{Milnor} gives, for invariant $X,Y\in V$,
\begin{equation}\label{eq:connection}
 \nabla_\nu\nu=0,\quad \nabla_\nu X=A_-X,\quad
 \nabla_X\nu=-SX,\quad \nabla_XY=\langle SX,Y\rangle\nu.
\end{equation}
These are also the connection formulas used in \cite{Moreno}. With $R(U,V)=[\nabla_U,\nabla_V]-\nabla_{[U,V]}$, the connection gives
\begin{align*}
 R(X,Y)Z&=-\langle SY,Z\rangle SX+\langle SX,Z\rangle SY,\\
 R(\nu,X)\nu&=(S^2+[S,A_-])X,\qquad
 R(\nu,X)Y=-\langle(S^2+[S,A_-])X,Y\rangle\nu.
\end{align*}
Taking traces gives $\operatorname{Ric}(X,Y)=-\langle[S,A_-]X,Y\rangle-(\tr S)\langle SX,Y\rangle$, zero mixed terms, and $\operatorname{Ric}(\nu,\nu)=-\tr(S^2)$. Since $\tr S=0$, this is \eqref{eq:ricci}. The torsion formula proves the final assertion directly.
\end{proof}

For a closed form, $\Delta_g\varphi=d\tau_2$, so \eqref{eq:torsion} also gives its initial Laplacian-flow velocity \cite{Karigiannis,LW,Lauret}. The evolving metric generally leaves the prescribed metric class.

When $S=0$, closedness gives $A=A_-\in\mathfrak{su}(V,J)$. This is the fixed-metric version of Freibert's parallel classification \cite[Theorem~3.6(a)]{Freibert}, also stated in \cite{Lauret}. The curvature formulas above give flatness directly, in agreement with \cite{Moreno}, where skew brackets give flat metrics even when the form is not closed. Freibert's more general treatment of torsion-free $H$-structures on almost-abelian solvmanifolds \cite{Freibert2026} places this parallel case in a broader setting. This agrees both with Milnor's flat metric Lie algebras \cite{Milnor} and with the theorem of Alekseevskii--Kimel'fel'd that Ricci-flat homogeneous Riemannian spaces are flat \cite{AK}.

If the metric of a nonempty fibre is Einstein, the trace on $V$ in \eqref{eq:ricci} makes the Einstein constant zero, because a commutator has trace zero. The $\nu\nu$ entry then forces $S=0$. Thus it is flat. Moreno's general almost-abelian Ricci formula \cite{Moreno} gives the same conclusion after imposing closedness. This is a specialization of the Einstein obstruction for calibrated $G_2$-structures on solvmanifolds proved by Fern\'andez--Fino--Manero \cite{FFM}.

\begin{proposition}\label{prop:dirac-aa}
Let $\sigma$ be the spin structure defined by the left-invariant orthonormal framing, on the simply connected group or on a compact quotient when one exists. For symmetric trace-free $A$, every spinor constant in this framing is harmonic for $D_\sigma$ and has constant length. If $A\neq0$ admits an invariant compatible closed form, the closed locus is nevertheless the proper subset of $\RP^7$ described in Corollary \ref{cor:symmetric}.
\end{proposition}
\begin{proof}
For $A_-=0$, lift \eqref{eq:connection} to spinors. If $s$ is constant in an invariant spin frame and $e_1,\ldots,e_6$ is an orthonormal basis of $V$, then
\[
 \nabla_\nu s=0,\qquad
 \nabla_{e_i}s=\tfrac12\sum_{j=1}^6 A_{ij}e_j\cdot\nu\cdot s,
 \qquad A_{ij}=\langle Ae_i,e_j\rangle.
\]
Consequently
\[
 D_\sigma s=\tfrac12\sum_{i,j=1}^6 A_{ij}e_i\cdot e_j\cdot\nu\cdot s
       =-\tfrac12(\tr A)\nu\cdot s=0.
\]
The off-diagonal terms cancel because $A$ is symmetric and Clifford multiplication anticommutes on orthogonal vectors. The unit invariant spinors modulo sign give all of $\RP^7$, whereas Corollary \ref{cor:symmetric} gives only one-dimensional or three-dimensional components for the closed locus when $A\neq0$. The harmonicity assertion also applies when $\Ccal_A$ is empty, for instance when a symmetric trace-free $A$ has an odd-dimensional eigenspace.
\end{proof}

On a compact quotient, Proposition \ref{prop:dirac-aa} concerns the spin structure defined by the descended invariant framing. Other spin structures, and spinors not constant in this framing, must be considered before identifying the full fibre. This distinction is already visible in flat quotients. The relation between flat solvmanifolds and compatible structures is studied by Tolcachier \cite{Tolcachier}.

Take $V=\R^2\oplus\R^4$ and $A=0\oplus2\pi I$, where $I$ is an orthogonal complex structure on $\R^4$. The standard semidirect-product coordinates identify $G_A$ isometrically with Euclidean $\R^7$, since $\exp(tA)$ is orthogonal. The invariant fibre is $\RP^3$ by row four of Table \ref{tab:normal}. For $\Gamma_1=\mathbb Z\ltimes\mathbb Z^6$, the action is by translations because $\exp A=\Id$. The quotient is a flat torus, whose full fibre is $\RP^7$ by \eqref{eq:all-parallel}. Hence the invariant fibre is a proper projective subspace of one component of the full fibre.

For $\Gamma_2=\tfrac12\mathbb Z\ltimes\mathbb Z^6$, the half-period acts by $\Id$ on $\R^2$ and by $-\Id$ on $\R^4$. Equivalently, the quotient is generated on $\R^7$ by translations in $\mathbb Z^7$ and
\[
 \gamma(x_1,\ldots,x_7)
 =(x_1+\tfrac12,x_2,x_3,-x_4,-x_5,-x_6,-x_7).
\]
It is compact and flat, with $b_1=3$. The holonomy lift
$\widetilde h=e_4e_5e_6e_7$ is a symmetric Clifford involution. It has trace zero, since conjugation by $e_4$ changes its sign, and therefore its two eigenspaces $E_+$ and $E_-$ each have dimension four. Translations must lift trivially for a parallel spinor to descend. The two choices $\pm\widetilde h$ then define precisely the two spin structures with parallel spinors. Equivalently, the projective lines fixed by the holonomy are exactly those in these eigenspaces. Hence,
\begin{equation}\label{eq:flat-two}
 \Ccal(g,o)=\mathbb P(E_+)\sqcup\mathbb P(E_-)
 \cong\RP^3\sqcup\RP^3.
\end{equation}
The invariant fibre is one of these components. Indeed, invariant closed forms are parallel and arise from the kernel of the infinitesimal spin action of $A$, which has dimension four. Their projective space is connected and has the dimension of a full component in \eqref{eq:flat-two}. The other component consists entirely of non-invariant forms. This also shows why the spin-sector qualification in \eqref{eq:lin} is necessary.

For symmetric brackets one can exclude non-invariant harmonic spinors after making the torus fibres sufficiently small. This gives a compact setting in which the algebraic model of Section \ref{sec:compact} is exactly the invariant model of this section.

\begin{theorem}\label{thm:collapsed}
Let $A=A^t$ with $\tr A=0$, and suppose $\exp(t_0A)$ preserves a lattice $\Lambda\subset V$ for some $t_0>0$. On
$M=(t_0\mathbb Z\ltimes\Lambda)\backslash G_A$, consider the metrics
\[
 g_\epsilon=dt^2+\epsilon^2\langle e^{-tA}dx,e^{-tA}dx\rangle.
\]
For all sufficiently small $\epsilon>0$, the Dirac kernel in the spin structure defined by the invariant framing consists exactly of the eight-dimensional space of invariant spinors. Every other spin structure has zero Dirac kernel. Consequently every compatible closed form is invariant and
\[
 \Ccal(g_\epsilon,o)\cong\Ccal_A.
\]
In particular, there are compact Riemannian seven-manifolds whose
entire compatible closed fibre is a disjoint union of four
projective circles, with every form in the fibre having nonzero
torsion.
\end{theorem}
\begin{proof}
Use the invariant orthonormal frame
\[
 \nu=\partial_t,\qquad
 E_a=\epsilon^{-1}(e^{tA}e_a)\cdot\partial_x,
 \qquad 1\leq a\leq6.
\]
The notation in the second expression means the directional derivative in the displayed coordinate vector. Proposition \ref{prop:dirac-aa} shows that the zero-order part of the Dirac operator in this frame is zero. On spinor-valued functions it is therefore
\[
 D=\gamma_\nu\partial_t+\sum_{a=1}^6\gamma_aE_a,
\]
where $\gamma_\nu,\gamma_a$ are constant Clifford matrices. Since $[E_a,E_b]=0$ and $[\partial_t,E_a]=\sum_b A_{ba}E_b$, its square is
\[
 D^2=-\partial_t^2-\sum_a E_a^2
       +\gamma_\nu\sum_{a,b}\gamma_a A_{ba}E_b.
\]
We now complexify the spinor bundle so that we can expand a
spinor in Fourier series along the torus fibres. Relative to
the global frame above, a spin structure is specified by a
sign character of the lattice group. The restriction of this
character to $\Lambda$ determines whether the spinor is
periodic or antiperiodic in each lattice direction. The
allowed frequencies therefore form a translate of
$\Lambda^*$ by a half-lattice vector, and in particular
every frequency belongs to $\tfrac12\Lambda^*$. For a frequency $k$, write
\[
 q(t,k)=2\pi\epsilon^{-1}e^{tA^t}k,
\]
using the metric to identify covectors with vectors. On that frequency the last term of $D^2$ is $i\gamma_\nu\gamma(Aq)$, whose operator norm is $|Aq|$. Thus the zeroth-order part is bounded below by
$|q|^2-\|A\|_{\mathrm{op}}|q|$.

There is a number $\delta>0$ such that
\[
 |e^{tA^t}k|\geq\delta
 \quad(0\leq t\leq t_0,\quad
       k\in\tfrac12\Lambda^*\setminus\{0\}).
\]
This follows from the positive minimum singular value on the compact interval and the discreteness of the half-lattice. Choose $\epsilon$ so that $2\pi\delta/\epsilon>\|A\|_{\mathrm{op}}$. Integration over a fundamental interval and Parseval's identity give, on the nonzero-frequency part $s_\perp$,
\[
 \|Ds_\perp\|_{L^2}^2
 \geq\|\partial_t s_\perp\|_{L^2}^2+
 \sum_{k\ne0}\int_0^{t_0}
 (|q(t,k)|^2-\|A\|_{\mathrm{op}}|q(t,k)|)|s_k(t)|^2\,dt.
\]
We have omitted the common positive volume factor from both
sides of the inequality. The identification at the endpoints
of the fundamental interval permutes the frequencies and
preserves the integrated norm. There is no boundary contribution,
since the monodromy preserves the lattice volume and the
identification of spinors is unitary. Thus the integration above
is precisely the identity
$\langle D^2s_\perp,s_\perp\rangle=\|Ds_\perp\|^2$ on the compact
quotient. By our choice of $\epsilon$, the lower bound is
strictly positive whenever $s_\perp\ne0$. A harmonic spinor
therefore has no nonzero-frequency part.

It remains to consider the zero-frequency part, which is present
only when the spin character restricts trivially to $\Lambda$.
For this part, the Dirac equation reduces to
$\gamma_\nu\partial_t s=0$, so $s$ is constant. Such a spinor
descends to the quotient only if the character of the base
generator is also trivial. This gives exactly the spin structure
defined by the invariant framing, whose Dirac kernel is the real
eight-dimensional space of constant spinors. By Theorem
\ref{thm:main}, every compatible closed form is therefore defined
by an invariant spinor satisfying the closure equation.
Rescaling all six fibre directions by the same factor $\epsilon$
leaves the orthonormal bracket matrix $A$ unchanged, so the full
closed fibre is $\Ccal_A$. The symmetric $2+2+2$ lattice
construction below proves the final assertion.
\end{proof}

The Fourier decomposition here is along the abelian fibres of a mapping torus. It serves the same purpose as representation-theoretic decompositions in calculations of harmonic spinors on nilmanifolds \cite{BMM,FGR}, but the positivity estimate above is specific to the symmetric almost-abelian metric. It does not assert that every metric on such a quotient has only invariant harmonic spinors.

\begin{remark}
Let $A=A^t$, $A\neq0$, and $\Ccal_A\ne\varnothing$. Then $V$ is the nilradical: a nilpotent ideal containing $c\nu+v$ with $c\ne0$ would contain $\im A$, on which its adjoint action is $cA$ and is not nilpotent. Bock's lattice results \cite{Bock} therefore imply that any lattice forces $\exp(t_0A)$ to preserve a full lattice in $V$ for some $t_0\ne0$. Conversely, this condition produces the cocompact subgroup $t_0\mathbb Z\ltimes\Lambda$ \cite{Bock}.

In the $4+2$ case the two eigenvalues of $A$ are $\lambda,-2\lambda$, with $\lambda\ne0$. Set $a=\exp(t_0\lambda)>0$. The characteristic polynomial of the integral lattice automorphism would be
\[
 P(x)=(x-a)^4(x-a^{-2})^2\in\mathbb Z[x].
\]
Since $a\ne a^{-2}$, it is the unique root of multiplicity four. Galois conjugation preserves multiplicities, so every conjugate of $a$ equals $a$, giving $a\in\mathbb Q$. Both $a$ and $a^{-1}$ are algebraic integers because the lattice automorphism is unimodular. Hence $a$ is an integer unit, and positivity forces $a=1$, a contradiction. Thus the nonzero symmetric $4+2$ groups admit no lattices.

The $2+2+2$ class does occur. Take any $M\in SL(3,\mathbb Z)$ with three distinct positive eigenvalues and choose an inner product making its eigenbasis orthonormal. Then
$A=(\log M)\oplus(\log M)$ is symmetric, trace-free, and has multiplicities $2+2+2$, while $\exp A=M\oplus M$ preserves $\mathbb Z^6$. Such matrices exist: the companion matrix of $x^3-5x^2+6x-1$ has determinant one and one root in each of $(0,1)$, $(1,2)$ and $(3,4)$. The resulting compact quotients have an invariant fixed-metric closed fibre consisting of four circles. After shrinking the torus fibres as in Theorem \ref{thm:collapsed}, these are the entire fibre. In the framing sector $\Hcal_\sigma=\R^8$ and $Z_\sigma$ is exactly the union of four projective lines described in Corollary \ref{cor:table}. All other $Z_\sigma$ are empty. Since $A\ne0$ is symmetric, Proposition \ref{prop:torsion} shows that every form in this fibre has nonzero torsion. 

Normal brackets with an invariant $\RP^3$ also occur on compact quotients with nonzero torsion. Let $M$ be the companion matrix of $x^3-x-1$. It has determinant one, one positive real eigenvalue $\mu>1$, and a non-real conjugate pair $\rho,\overline\rho$. Choose a real logarithm and an inner product making $\log M$ normal, and set $A=\log M\oplus\log M$. Its spectrum is row four of Table \ref{tab:normal}, because $\log\mu=-2\log|\rho|$. The real part is nonzero and $\exp A=M\oplus M$ preserves a lattice. Thus its invariant closed fibre is $\RP^3$ with nonzero torsion. No assertion about the full non-invariant fibre is made for this metric. Compact solvable calibrated examples and lattice constructions also appear in Fern\'andez \cite{FernandezCompact} and Manero \cite{Manero}.
\end{remark}

\printbibliography

\end{document}